\input amstex
\advance\vsize-0.5cm\voffset=-0.5cm\advance\hsize1cm\hoffset0cm
\magnification=\magstep1
\documentstyle{amsppt}
\NoBlackBoxes
\topmatter
\title
{Surgery spectral sequence and stratified manifolds}
\endtitle
\author Rolando Jimenez, Yurij V. Muranov,  and  Du\v san Repov\v s
\endauthor
\leftheadtext{R. Jimenez, Y.V. Muranov, and D. Repov\v s} 
\thanks{First author was partially supported by   CONACyT, DGAPA-UNAM,
Fulbright-Garcia Robles and UW-Madi\-son grant, second author by 
Russian Foundation for Fundamental
Research Grant no. 02--01--00014, and the third author by the 
Ministry of Education,
Science and Sport of the Republic of Slovenia research programme 
No. P1--0292--0101--04.}
\endthanks
\endkeywords
\abstract{Cappell and Shaneson pointed out in 1978 
interesting properties of 
Browder - Livesay invariants which are similar to differentials in some
spectral sequence. Such spectral sequence was constructed in 1991 by
Hambleton and Kharshiladze. This spectral sequence is  closely related 
to a problem of realization of elements of Wall groups by  normal
maps of closed manifolds.  The main step of construction 
of the spectral sequence is an infinite filtration of spectra
in which only the first two, as is well-known, 
have clear  geometric  sense.
The first one is a spectrum $\Bbb L(\pi_1(X))$ 
for surgery obstruction groups of a manifold $X$ and the second  
$\Bbb LP_*(F)$ is a spectrum for surgery on a Browder-Livesay 
manifold pair
$Y\subset X$. The geometric sense of the third term of filtration was
explained by Muranov, Repov\v s, and Spaggiari in 2002. 
In the present paper we give a geometric interpretation of all spectra
of filtration in construction of Hambleton and Kharshiladze.  
We introduce groups of obstructions to surgery on a system of embeddedd 
manifolds and prove that  spectra which  realize  
these groups coincide with spectra in the  filtration of Hambleton and
Kharshiladze. 
We describe algebraic and geometric properties of
introduced obstruction groups and their  relations to the classical 
surgery theory.  We prove isomorphism between introduced groups
and Browder-Quinn $L$-groups of stratified manifolds.   
We give an application 
of our results to closed manifold surgery problem and iterated 
Browder-Livesay invariant.}
\endabstract
\endtopmatter
\document

\bigskip
\subhead 1. Introduction
\endsubhead
\bigskip

Let $X^n$ be a  closed $n$@-dimensional $CAT$ ($CAT=TOP, PL,  DIFF$)
manifold with the fundamental group $\pi =\pi_1(X)$ which is given 
with  a homomorphism 
of the orientation $w:\pi_1(X)\to \{\pm 1\}$. In the sequel  
we shall assume that 
all groups are given with an orientation homomorphism 
and shall not specify this in notations without necessarity. 

A fundamental problem of geometric topology is
 to describe all
closed $n$@-dimensional $CAT$@-manifolds which are homotopy
(simple homotopy) equivalent to $X$. More precisely, 
let $f:M^n\to X^n$ be
a simple homotopy equivalence of $CAT$-manifolds. 
The structure set  $\Cal S^{CAT}(X)$ is the set of
$s$@-cobordism classes of equivalence of $CAT$@-manifolds which
are simple homotopy equivalent to  $X^n$ (see \cite{34},\cite{29} and
\cite{30, p. 542}). 
The elements of $\Cal S^{CAT}(X)$  are called 
$s$@-triangulations of the manifold $X$.

 The Sullivan-Novikov-Wall surgery exact sequence
 $$
\cdots \to L_{n+1}(\pi) \rightarrow
\Cal S^{CAT}(X) \to [X, G/CAT]  \overset{\sigma}\to{\rightarrow} L_n(\pi)\cdots
\tag 1.1
$$
is the main tool
for describing  the structure set $\Cal S^{CAT}(X)$ (see \cite{34} and
\cite{30}).

Hereafter we shall  consider only topological
manifolds ($CAT=TOP$)  and  groups $L_*(\pi)=L^s_*(\pi)$ which
give obstructions to simple homotopy equivalence
(see \cite{34, \S 10} and \cite{30}). 
To describe the structure set $\Cal S^{TOP}(X)$ we must compute
the set of normal  invariants $[X, G/TOP]$, the surgery obstruction 
groups $L_{n}(\pi)$ and the map $\sigma$ in (1.1).   
To describe  the map $\sigma$ we must know 
what elements of the group $L_n(\pi)$ are realized 
by  normal maps of closed manifolds.  

The algebraic surgery exact sequence of Ranicki 
(see \cite{29} and \cite{30}) 
$$ 
\cdots\rightarrow L_{m+1}(\pi_1(X))\rightarrow \Cal 
S_{m+1}(X)\rightarrow 
H_m(X; \bold L_{\bullet})\overset{\sigma}\to{\rightarrow} 
 L_m(\pi_1(X))\rightarrow \cdots 
\tag 1.2 
$$ 
is defined for any topological space $X$.  In particular, it defines 
an assembly map 
$$
H_n(K(\pi, 1); \bold L_{\bullet})\overset{A}\to{\rightarrow} 
 L_n(\pi) 
\tag 1.3 
$$
and $Image (A)\subset L_n(\pi)$ is the subgroup consisting of the elements 
which can be realized by normal maps of closed manifolds
(see, for example, \cite{34, \S 13}).

 If the space $X$ is 
simple homotopy equivalent to a topological $n$-dimensional 
manifold $n\geq 5$, 
then the  exact sequence (1.1) is isomorphic to corresponding part of 
(1.2). Exact sequence (1.2) is realized on the spectra 
level by a map of spectra     
$$ 
X_+\land \bold L_{\bullet}\rightarrow \Bbb L(\pi_1(X)) 
\tag 1.4 
$$ 
where $\Bbb L(\pi_1(X))$ is the surgery $L$-spectrum of 
 the fundamental group $\pi_1(X)$ with 
$$ 
\pi_n(\Bbb L(\pi_1(X)))\cong L_n(\pi_1(X)) 
$$ 
 and 
$\bold L_{\bullet}$ is the 1-connected cover of the 
$\Omega$-spectrum 
 $\Bbb L(\Bbb Z)$ 
 such that $\bold {L_{\bullet}}_0 \simeq G/TOP$. 
 
In particular, for  manifold $X$ we have 
isomorphisms 
 $\Cal S_{n + 1} (X)\cong \Cal S^{TOP}(X)$ and  
$H_n(X; \bold L_{\bullet}) \cong [X, G/TOP]$ 
(see \cite{29} and \cite{30}).

Approaches to computation of the structure set 
$\Cal S^{TOP}(X)$  are very different   for
the cases finite and infinite group $\pi$. The case 
of infinite group is closely related to the Novikov 
Conjecture (see, for example \cite{9}). In the case 
of finite groups the solution of the problem 
for the special case of decorations (the case 
of intermediate groups $L^{\prime}$) was given in \cite{12}.
The fundamental 
 results of \cite{12} are based on analysis of assembly
map and methods of \cite{5} and \cite{10}.  
The methods developed in  \cite{4}, \cite{5}, \cite{10}, 
and \cite{16} make it possible  to prove the nonrealizability 
of elements of the Wall group $L_n(\pi)$ which do not lie 
in the image of the natural map $L_n(1)\to L_n(\pi)$ for arbitrary
case. In particular,
Hambleton solved in \cite{10} the corresponding 
problem for projective Novikov groups $L^p_*$.    
These  methods are mostly algebraic and are 
based on the algebraic theory 
of 
splitting homotopy equivalence along submanifolds.  
          
Let  $Y\subset X$ be a submanifold of codimension $q$ in a closed 
topological  manifold $X$ of dimension $n$.
A simple homotopy equivalence $f: M\to X$ splits along the submanifold $Y$
if it is
homotopy equivalent to  a map $g$ which is transversal to $Y$ with  
a submanifold $N=g^{-1}(Y)\subset M$ and the restrictions
$$
g|_N: N\to Y,\ \  g|_{(M \setminus N)}:M\setminus N\to X\setminus Y
\tag 1.5
$$
 are simple homotopy equivalences.
Let  $U$ be a tubular neighborhood of
the submanifold $Y$ in $X$ with boundary $\partial U$. Denote
by 
$$
F=
\left(\matrix
\pi_1(\partial U)& \to &\pi_1(X\setminus Y)\\
\downarrow  & &     \downarrow\\
\pi_1(U)&\to &\pi_1(X)\\
\endmatrix\right)
\tag 1.6 
$$
a push-out square of fundamental groups 
with  orientations.
There exists a group $LS_{n-q}(F)$ of
obstructions to splitting (see \cite{34} and \cite{30}) which
depends only on $n-q \bmod 4 $ and the square $F$.

Let  $(f, b) : M\to X$ be a normal map with 
 $b: \nu_M\to \xi$  a map of fibrations 
covering $f$ where $\xi$ is  a topological reduction 
of the Spivak normal fibration over $X$ (\cite{29} and 
\cite{30}). In this case an obstruction for existence 
in the normal bordism class of the map $(f, b)$
a map $(g, b^{\prime})$ with properties (1.5) lies in
the group 
$LP_{n-q}(F)$ of obstructions to surgery on  manifold pairs
(see \cite{34} and \cite{30}). Also this group
depends  only on $n-q\bmod 4$ and the square $F$ of fundamental 
groups. 

The main relation between $LS_*$- and $LP_*$-groups  
and algebraic surgery exact
sequence (1.2) is given by the following commutative diagram
\cite{34, \S 11}  
$$ 
\matrix 
\dots \rightarrow &\Cal S_{n+1}(X)&\rightarrow 
&H_n(X; \bold L_{\bullet})&\overset{\sigma}\to{\rightarrow} 
& L_{n}(\pi_1(X))&\to\cdots \\ 
&\downarrow & 
& \downarrow \sigma_1 & &\downarrow = & \\ 
\dots \rightarrow &LS_{n-q}(F) &\rightarrow 
&LP_{n-q}(F)&\overset{s}\to{\rightarrow} 
&L_{n}(\pi_1(X))&\overset{\partial}\to{\to}\cdots \\ 
\endmatrix 
\tag 1.7 
$$ 
where the  rows are exact sequences. 
It follows from (1.7) that the image of the map $\sigma$ 
lies in the kernel of the map
$$
\partial:  L_{n}(\pi_1(X)) \to LS_{n-q-1}(F)
$$
which has a clear  geometric description (\cite{5} and \cite{16}).

The bottom row of diagram (1.7) 
fits the following 
braid of exact sequences (see \cite{34, p. 264} and \cite{30, \S 7.2}) 
$$ 
 \matrix \rightarrow & {L}_{n+1}(C) & \longrightarrow & 
   L_{n+1}(D) & 
\overset{\partial}\to{\rightarrow} & LS_{n-q}(F)&\rightarrow 
 \cr 
\ & \nearrow \  \  \ \ \ \searrow &\ &^s\nearrow \ \ \ \ \ 
                           \searrow 
   & \ & \nearrow \ \ \  \ \ \searrow & \ \cr 
  \ & \ & LP_{n-q+1}(F)& \ & L_{n+1}(C\to D) & \ & \ \cr 
   \ & \searrow \ \ \ \ \ \nearrow &\ &\searrow \ \ \ \ \ 
\nearrow 
 & \ & \searrow \ \    \ \ \ \nearrow & \ \cr 
\rightarrow & LS_{n-q+1}(F) & \longrightarrow & 
 L_{n-q+1}(B) & 
 \longrightarrow & {L}_{n}(C) & \rightarrow 
\endmatrix 
\tag 1.8 
 $$ 
where $A=\pi_1(\partial U)$, $B=\pi_1(Y)$, 
$C=\pi_1(X\setminus Y)$, and $D=\pi_1(X)$. 
                  
Now let a pair of manifolds $(X,Y)$ be a Browder-Livesay pair (\cite{2},
 \cite{5},
 \cite{10}, \cite{16}, and \cite{22}). 
This means that $Y$ is a one-sided submanifold 
of codimension  1
of the  manifold $X$ and the natural embedding  $Y \to X$ induces 
an isomorphism of the fundamental groups. 
In this case the square $F$  of fundamental groups (1.6) 
has the  following form 
$$ 
F= 
\left(\matrix 
\pi_1(\partial U)& \to &\pi_1(X\setminus Y)\\ 
\downarrow & & \downarrow\\ 
\pi_1(Y)&\to &\pi_1(X)\\ 
\endmatrix\right)= 
\left(\matrix 
A\ \ \ &\overset{\cong}\to{\to} &A\ \ \ \\ 
\downarrow i_- & & \downarrow i_+\\ 
B^-&\overset{\cong}\to{\to} &B^+\\ 
\endmatrix\right). 
\tag 1.9
$$ 
The  orientation  of the group $B^-$ in (1.9) 
 differs from the orientation of the group 
$B^+$ outside the images of the vertical maps (which are 
inclusions of index 2).  All maps in the square (1.9), 
except the lower horizontal  map, 
preserve the orientation. 
The lower isomorphism preserves the orientation on the image 
of $i_-$, 
and reverses the orientation outside this image. 
In this case, we have an  
isomorphism  
$$ 
LP_n(F)\cong L_{n+1}(i^*_-), 
$$ 
where $i^*: L_{n+1}(B^-)\to L_{n+1}(A)$ is the transfer map. 
The group  
$LS_*(F)$ is denoted by $LN_*(A\to B^+)=LN_*(A\to B)$ \cite{34}, and 
is called the Browder-Livesay group. 

Cappell and Shaneson proved \cite{9} that for a Browder-Livesay
pair $(X,Y)$ the elements which do not lie in the kernel  
of the map 
$$
\partial:  L_{n}(\pi_1(X)) \to LN_{n-2}(\pi_1(X\setminus Y)\to \pi_1(X))
$$
cannot be realized by a normal map of closed manifolds. 

The diagram (1.8) for Browder-Livesay pairs has 
an algebraic description (see \cite{10} and \cite{31}).
This diagram was investigated from 
algebraic and geometric point of view 
in several papers (see 
\cite{10}, \cite{12}, \cite{13}, \cite{16}, \cite{20}, \cite{21}, 
 \cite{22},\cite{23}, \cite{24}, 
 and \cite{31}).  

Subsequently  a spectral sequence in surgery theory was 
 constructed in \cite{13} using realization of commutative diagram (1.8)
for a Browder-Livesay pair 
on the spectra level. Consider the filtration of spectra from \cite{13} 
$$ 
\cdots\rightarrow \Bbb X_{3,0}\rightarrow \Bbb X_{2,0} 
\rightarrow \Bbb X_{1,0}\rightarrow \Bbb X_{0,0}\rightarrow 
\Bbb X_{-1,0}\rightarrow \cdots 
\tag 1.10  
$$ 
where $\Bbb X_{0,0}=\Bbb L(\pi_1(X))$ is a surgery spectrum  
with $\pi_n(\Bbb L(\pi_1(X)))= L_n(\pi_1(X))$
and  $\Bbb X_{1,0}$ is a spectrum for surgery obstruction 
groups on the manifold pair $(X, Y)$
$$
\Sigma\Bbb LP(F)=\Bbb L(i^*_-).
$$ 

The map $s$ in commutative diagrams (1.7) and (1.8) is induced by the
map of spectra 
$\Bbb X_{1,0}\rightarrow \Bbb X_{0,0}$ from filtration (1.10).
Another spectra of filtration is defined inductively using the 
pullback
construction and as is well-known they have no geometric meaning. 
It follows from \cite{13} that the surgery spectral sequence 
is closely   related to the iterated Browder-Livesay invariants 
and to the oozing problem.
Other versions of surgery spectral exact sequence were 
obtained in papers   \cite{6}, \cite{7}, 
\cite{15}, and  \cite{20}.

Let $Z\subset Y \subset X$ be a triple of closed topological 
manifolds such that $n$ is the dimension of $X$, $q$ is the 
codimension of $Y$ in $X$, and $q'$ is the  codimension 
of $Z$ in $Y$.   Groups of 
obstructions to surgery $LT_{n-q-q^{\prime}}(X,Y,Z)$ 
on  manifold triples 
were introduced in  \cite{26}. These groups are
realized on the spectra level by a spectrum 
$\Bbb LT(X,Y,Z)$
and they are a natural 
generalization of surgery obstruction groups $LP_*$
for manifold pairs. The natural forgetful map 
$t:LT_{n}(X,Y,Z)\to LP_{n+1}$ 
which is realized on the spectra level is well-defined.

If the triple $(X,Y,Z)$ consists of the  Browder-Livesay's pairs $(X,Y)$ 
and $(Y,Z)$,  then  the spectrum   $\Sigma^2\Bbb LT(X,Y,Z)$ 
coincides with the  spectrum $\Bbb X_{2,0}$ of filtration 
(1.10).  The map $\Bbb X_{2,0}\to \Bbb X_{1,0}$ of filtration (1.10) coincides 
with the map $t$ on the spectra level \cite{26}.   

Now let 
$$
X_k\subset X_{k-1}\subset \cdots \subset X_2\subset X_1\subset X_0=X
\tag 1.11
$$
be a filtration $\Cal X$ of a closed topological manifold $X$ by
locally flat embedded submanifolds. 
Denote by $l_j$ the dimension 
of the submanifold $X_j$ and by $q_j$ the codimension of $X_j$ 
in $X_{j-1}$ for $1\leq j\leq k$. 
 We shall suppose that 
every pair of manifolds from (1.11) is a topological manifold
pair in the sense of Ranicki \cite{30, \S 7.2} and that dimension 
$l_k\geq 5$. 

For every nonempty subset $B\subset \{k,k-1, \dots , 2,1,0\}$ 
filtration (1.11) defines the restricted filtration $\Cal X_B$
which is obtained by forgetting the submanifolds $X_j$ from 
filtration (1.10) with $j\in \{k,k-1, \dots , 2,1,0\}\setminus B$.  
In particular a restricted filtration    
$$
X_j\subset X_{j-1}\subset \cdots \subset X_2\subset X_1\subset X_0=X
\tag 1.12
$$
where $B=\{j,j-1, \dots , 2,1,0\}$ 
for every 
$0\leq j\leq k$ is well-defined. We shall denote  the restricted 
filtration (1.12) by $\Cal X_j$.
 
 For a  
simple homotopy eqiuvalence $f:M\to X$ we define a concept 
of  an $s$-triangulation 
of the filtration (1.11) 
in section 2
and prove several technical results.  
In particular,  we prove that for the manifold triple 
$Z\subset Y\subset X$ the surgery obstruction groups 
$LT_*(X, Y, Z)$ from \cite{26} coincide with the 
Browder-Quinn groups $L^{BQ}$ (see \cite{3} and \cite{35}) 
of the stratified manifold
$Z\subset Y\subset X$.    

We then construct in Section 3 the groups of 
obstructions to $s$-triangulation of a filtration $\Cal X$ (1.11)
of embedded manifolds
and study their  properties.   
We introduce  obstruction groups $LM^j_{i}(\Cal X) 
\ (0\leq j\leq k)$
which have period $4$ for subscript $i$ and which
are realized on the spectra level by spectra $\Bbb LM^j(\Cal X)$ 
with $\pi_i(\Bbb LM^j(\Cal X))=LM^j_i(\Cal X)$. 
The groups $LM^k_*(\Cal X)$ coincide with Browder-Quinn stratified 
$L$-groups $L^{BQ}_*(\Cal X)$
(see \cite{3} and  \cite{35})
 up to a shift of dimension $*$.  
The spectrum $\Bbb LM^0$
coincides with the spectrum $\Bbb L(\pi_1(X))$, 
the spectrum $\Bbb LM^1$
coincides with the spectrum $\Bbb LP(F)$ for the pair $(X_0, X_1)$
(see \cite{30} and \cite{34}),
and the spectrum $\Bbb LM^2$
coincides with the spectrum $\Bbb LT$ for the triple 
$(X_0, X_1, X_2)$ (see \cite{26} and \cite{28}).

Let $(f,b): (M\to X) $ be  a normal map to the manifold $X$ 
with the filtration (1.11). For groups introduced above and 
$0 \leq j\leq k $  an obstruction 
$\Theta_{j}(f)\in  LM^j_{l_j}$ is defined. It is proved in Theorem 3.9 that
this obstruction   
is trivial if and only if the map $f$ is normally bordant 
to an  $s$-triangulation of the restricted filtration $\Cal X_j$ (1.12).

In section 3 we define the  natural forgetful maps  
$$
LM^k_{l_k}\to LM^{k-1}_{l_{k-1}}\to \cdots
\to LM^{1}_{l_1}\to LM^{0}_{l_{0}}
\tag 1.13
$$
which are realized on the spectra level
by 
maps of spectra 
$$
\Sigma^{n-l_k} \Bbb LM^k\to
\Sigma^{n-l_{k-1}} \Bbb LM^{k-1}\to \cdots
\to\Sigma^{n-l_1} \Bbb LM^1 \to \Bbb LM^{0}.
\tag 1.14
$$

Ranicki introduced in \cite{30} 
a set $\Cal S_{n+1}(X, Y, \xi)$ of
homotopy triangulations  of a pair of manifolds $(X, Y)$, where
$\xi$ denotes the normal bundle of $Y$ in $X$. This set consists
of concordance classes of maps  $f:(M, N)\to (X,Y)$ which are
split along $Y$. This structure set is a natural generalization 
of the structure set $\Cal S_{n+1}(X)$  from exact sequence (1.2) 
and fits into the   exact
sequence
 (see \cite{30, \S 7.2})
$$
\cdots \to  \Cal S_{n+1}(X,Y,\xi)  \to
H_{n}(X, \bold L_{\bullet}) \to LP_{n-q}(F)\to\cdots 
\tag 1.15 
$$
which is a natural generalization of (1.2) to the case of
manifold pairs. 

In Section 3  we introduce structure sets   for the filtration (1.11) 
which generalize 
structure sets 
$\Cal S_{n+1}(X, Y, \xi)$ and  $\Cal S_{n+1}(X)$  
and we study their properties. Some results
for the case of manifold triples were obtained 
in \cite{26},  \cite{27}, and \cite{28}.  
\bigskip

Let  all pairs $X_{i+1}\subset X_{i}$ 
in (1.11) be  Browder-Livesay
pairs for $0\leq i\leq k-1$. In Section 4 we apply our results to 
an investigation of iterated Browder-Livesay invariants and we describe 
relations of introduced groups to
surgery spectral sequence.
It is proved in Theorem 4.1 that in this case 
filtration (1.14) coincides with the left 
part starting with $\Bbb X_{0,0}$ of filtration (1.10) 
for spectral sequence of Hambleton and 
Kharshiladze. Furthemore in Section 4  we investigate relations 
of  groups $LM^i_*$ to the realization of
elements of Wall groups by normal maps of closed manifolds.

\subhead  2. Preliminaries and technical results
\endsubhead
\bigskip

In this section we recall some preliminary results about 
surgery on topological manifolds and use of  surgey $L$-spectra
(see \cite{1}, \cite{8}, \cite{11},
 \cite{26}, \cite{29}, \cite{30}, and \cite{33}). 
We shall give the necesary definitions and 
prove several technical results.

We shall consider a case of topological manifolds  and follow
notations from \cite{30, \S 7.2}. 
Let $(X, Y, \xi)$ be a codimension $q$   manifold 
pair  in the sense of Ranicki (see \cite{30, \S 7.2}), 
i. e.  
a locally flat closed 
submanifold  $Y\subset X$  given with a normal fibration 
 $$ 
  \xi=\xi_{Y\subset X}: Y \to \widetilde{BTOP}(q) 
 $$ 
with the associated $(D^q, S^{q-1})$ 
fibration 
$$
(D^q, S^{q-1})\to (E(\xi), S(\xi))\to Y
\tag 2.1
$$ 
and we have a decomposition of the closed manifold  
$$ 
 X =E(\xi)\cup_{S(\xi)}\overline{X\setminus E(\xi)}. 
$$ 
A topological normal map  \cite{30, \S 7.2}
$$ 
((f,b), (g,c)):(M,N)\to (X,Y) 
$$
to the manifold pair $(X, Y, \xi)$) is 
represented by a normal map 
$(f,b)$ to the manifold $X$ 
which is transversal to $Y$ with $N=f^{-1}(Y)$, and 
$(M,N)$ is a topological manifold pair with a normal 
fibration 
$$ 
\nu:N\overset{f|_N}\to{\to} Y \overset{\xi}\to{\to} 
\widetilde{BTOP}(q). 
$$ 
Additionaly, the following conditions are satisfied:

(i) the restriction 
$$ 
(f,b)|_N =(g,c) :N\to Y 
$$ 
is a normal map;

(ii) the restriction 
$$ 
(f,b)|_P =(h,d) :(P, S(\nu))\to (Z,S(\xi)) 
$$ 
is a normal map to the pair $(Z,S(\xi))$, 
where 
$$ 
P=\overline{M\setminus E(\nu)}, \ \ Z=\overline{X\setminus 
E(\xi)}; 
$$ 

(iii) the restriction 
$$ 
(h,d)|_{S(\nu)}: S(\nu)\to S(\xi) 
$$ 
coincides with the induced map 
$$ 
(g,c)^!: S(\nu)\to S(\xi), 
$$ 
and $(f,b)=(g,c)^!\cup(h,d)$. 

The normal maps to $(X,Y,\xi)$ are called $t$-triangulations 
of the manifold pair $(X,Y)$ 
and the set of concordance classes of $t$-triangulations 
of the pair $(X, Y, \xi)$ coincides with the set of $t$-triangulations 
of the manifold $X$ \cite{30, Proposition 7.2.3}. 

 An $s$@-triangulation of a manifold pair $(X,Y, \xi)$  in 
topological category 
 \cite{30, p. 571}
is a $t$@-triangulation 
 of this pair for which the maps 
 $$ 
 f:M\to X,\ g:N\to Y, \ 
 \text{and} \ 
 (P, S(\nu))\to (Z, S(\xi)) 
\tag 2.2
 $$ 
are simple homotopy equivalences ($s$@-triangulations). 

A simple homotopy equivalence $f: M\to X$ splits along a 
submanifold   $Y$ if it is homotopy equivalent to a map $g$ which is 
$s$-triangulation of $(X, Y, \xi)$ i.e. it satisfies conditions
(2.2). In this case $f$ 
represents an element of  $\Cal S_{n+1}(X, Y, \xi)$.
It follows from the definition of $s$@-triangulation of the pair
$(X,Y, \xi)$ that the forgetful maps
 $$
 \matrix
 \Cal S_{n+1}(X, Y, \xi) \to \Cal S_{n+1}(X), & (f, g)\to f; \cr
 \ & \cr
 \Cal S_{n+1}(X, Y, \xi) \to \Cal S_{n-q+1}(Y), &(f,g)\to g \cr
 \endmatrix
$$
are well-defined. In the general case the map $ \Cal S_{n+1}(X, Y,
\xi) \to \Cal S_{n+1}(X) $ is not an epimorphism or a monomorphism
\cite{30, p. 571}.

Consider a triple $Z^{n-q-q^{\prime}}\subset Y^{n-q} \subset X^n$
of closed topological manifolds.   We shall assume that every
submanifold is locally flat in the ambient manifold and that it is
equipped by  the structure of the normal topological bundle (see
\cite{30, pages 562--563} and \cite{26}). 
Every pair of manifolds defines the following topological normal bundles
which we denote in  the following way:
  $\xi$  for the submanifold $Y$ in $X$,
  $\eta$ for the submanifold  $Z$ in  $Y$, and
  $\nu$ for the submanifold $Z$ in $X$.
 We denote
the spaces with boundaries of associated fibrations (2.1) 
by $(E(\xi),S(\xi))$, $(E(\eta), S(\eta))$, and $(E(\nu), S(\nu))$,
respectively. 
Let $\xi|_{E(\eta)}$ be a  
restriction of the bundle $\xi$ on a space
$E(\eta)$ of normal bundle $\eta$ 
with a   restriction of fibration (2.1)  
$$
(D^q, S^{q-1})\to (E^{\prime}(\xi), S^{\prime}(\xi))\to E(\eta)
$$
and 
$\xi|_{S(\eta)}$ be a restriction with a restriction 
of fibration 
$$
(D^q, S^{q-1})\to (E^{\prime\prime}(\xi), S^{\prime\prime}(\xi))\to S(\eta)
$$
We assume that the space $E(\nu)$ of the 
normal bundle  $\nu$ is identified with the space $E^{\prime}(\xi)$
of the restriction
 $\xi|_{E(\eta)}$ in such a way that the following condition 
on the boundary is satisfied
$$
S(\nu)=E^{\prime\prime}(\xi)\cup S^{\prime}(\xi).
\tag 2.3
$$

\proclaim{Remark 2.1} The existence of normal bundles of the submanifolds 
for the manifold triple $Z^{n-q-q^{\prime}}\subset Y^{n-q} \subset X^n$
with the associated fibrations with  conditions (2.3) 
implies that the triple
$Z\subset Y\subset X$   is a $\Cal C$-stratified set in the 
sense of Browder and Quinn \cite{3}. 
\endproclaim 
\smallskip

Denote by $\Cal X$ a filtration of a closed manifold $X^n$ by a 
system 
of submanifolds (1.11). All pairs are given together 
with  normal bundles and corresponding $(D^*, S^{*-1})$ fibrations (2.1).
We shall suppose that for every triple of manifolds 
$X_j\subset X_l\subset X_m$ with $k\geq j> l> m\geq 0$
the conditions on the normal bundles similar to (2.3) for
the triple $Z\subset Y\subset X$ are satisfied.    

\proclaim{Remark 2.2} Under the assumptions above the 
filtration (1.11) $\Cal X$
gives a $\Cal C$-stratified set in the 
sense of Browder and Quinn \cite{3} --- this  follows from Remark 2.1
and Definition 4.2 from \cite{3}. 
\endproclaim

A codimension $q$  
manifold pair with boundaries 
$(Y, \partial Y)\subset (X, \partial X)$  
is defined in  \cite{30, p. 585}.
We  have a normal fibration 
$(\xi, \partial \xi)$ over the pair $(Y, \partial Y)$
and a decomposition  
$$
(X, \partial X)= (E(\xi)\cup_{S(\xi)}Z, 
E(\partial\xi)\cup_{S(\partial\xi)}\partial_+Z)
\tag 2.4
$$
where $(Z; \partial_+Z, S(\xi); S(\partial\xi))$ is a manifold
triad. Note  that here $\partial_+Z =\overline{\partial X\setminus 
E(\partial\xi)}$.

A topological 
normal map  
of manifold pairs with boundaries 
$$
(f, \partial f): (M, \partial M)\to (X, \partial X)
$$
provides  
a normal fibration 
$(\nu, \partial \nu)$ over the pair $(N, \partial N)$ 
(see \cite{30, p. 570})
where 
$$
(N, \partial N)=(f^{-1}(Y),(\partial f)^{-1}(\partial Y)).
$$ 
We have the 
following decomposition 
$$
(M, \partial M)= (E(\nu)\cup_{S(\nu)}P, 
E(\partial\nu)\cup_{S(\partial\nu)}\partial_+P)
\tag 2.5
$$
where $(P; \partial_+P, S(\nu); S(\partial\nu))$ is a manifold
triad.  

Now we define filtration $(\Cal X, \partial\Cal X)$ for the case 
of manifolds with boundaries as filtration  
$$
(X_k, \partial X_k)\subset 
(X_{k-1},\partial X_{k-1})\subset 
\cdots \subset (X_0, \partial X_0)=(X, \partial X)
\tag 2.6
$$
where all constituent pairs of manifolds with boundaries  satisfy  
properties which are similar to (2.4).
We also assume  that normal bundles of manifolds of filtration 
and of boundaries satisfy  properties 
similar to (2.3). 

\proclaim{Remark 2.3} Under the assumptions above the 
filtration (1.11) $\Cal X$
yields  a filtration of manifolds with boundaries
$$
\left(X_{k-1}\setminus X_k, \partial(X_{k-1}\setminus X_k)\right)  \subset 
\left(X_{k-2}\setminus X_{k},\partial(X_{k-2}\setminus X_k)\right)\subset 
\cdots \subset \left(X \setminus  X_k,\partial(X\setminus X_k)\right)
$$
This filtration is a
$\Cal C$-stratified manifold with boundary in the 
sense of  \cite{3} and  \cite{35}. We shall denote this filtration by 
$\overline{\Cal X_k}=\overline{\Cal X}$. In a similar way we can 
construct a filtration 
$\overline{\Cal X_j}$ using restricted filtration (1.12).
\endproclaim

\proclaim{Definition} A topological normal map to the filtration $\Cal X$ 
(1.11)
($t$-triangulation of the filtration $\Cal X$)
is a topological normal map 
$(f, b):M\to X$ which is topologically transversal to 
every submanifold of filtration with transversal preimages 
$M_0=M, M_i=f^{-1}(X_i)$ for $0\leq i\leq k$. We shall 
additionaly assume that  restriction 
on every pair of submanifolds $(M_j, M_l), \ (j\geq l)$ is topological normal 
map to the manifold pair $(M_j, M_l)$. 
In a natural way we can define the bordism of such 
maps, and the bordism classes are denoted by $\Cal T(\Cal X)$
(see \cite{3} and \cite{35}). 
\endproclaim

It is clear that
a $t$-triangulation of the filtration $\Cal X$ gives 
a $t$-triangulation of  
a restricted filtration $\Cal X_B$ for 
every nonempty subset $B\subset \{k,k-1, \dots , 2,1,0\}$.
In particular for every submanifold $X_j$ from the given filtration 
we have a forgetful map 
of $\Cal T(\Cal X)$ to the set $[X_j, G/TOP]$ of normal maps to 
the manifold
$X_j$. 
 
\proclaim{Proposition 2.4} (\cite{3} and \cite{30}) The natural forgetful map 
$\Cal T(\Cal X)\to [X, G/TOP]$ is an isomorphism. 
\endproclaim

\demo{Proof} Topological transversality (see 
\cite{3},  \cite{30, Proposition 7.2.3}, and \cite{35}) 
and induction on the number of elements of filtration. \qed
\enddemo
\smallskip

\proclaim{Definition} A $t$-triangulation 
$
(f,b): M\to X
$
of the filtration $\Cal X$ (1.11)
is an  $s$@-triangulation of the filtration   $\Cal X$ if the
constituent normal maps of pairs 
$$
(M_j, M_l)\to  (X_j, X_l),  0\leq j< l\leq k
$$
are
$s$@-triangulations i.e. they satisfy  the properties 
which are similar to properties (2.2) for the manifold pair $(X,Y)$.
\endproclaim  
\smallskip

\proclaim{Proposition 2.5}
Let  $t$-triangulation 
$(f,b):\Cal M\to \Cal X$  define an  $s$-triangulation 
$f_k: \Cal M_k\to \Cal X_k$ where filtration
$\Cal X_k$ is obtained from $\Cal X$ by forgetting 
the submanifold $X_k$ and similarly for the
$\Cal M_k$. Suppose that the
restriction $f|_{M_k}$ is an $s$-triangulation of the pair 
$(X_{k-1}, X_k)$.   Then $(f,b)$ is an $s$-triangulation 
of $\Cal X$.
\endproclaim

\demo{Proof} It is suffices to prove that for every 
submanifold $X_j\subset X, \ 0\leq j\leq k-2$ the restricted map 
$$
f|_{M_j\setminus M_k}: (M_j\setminus M_k)\to X_j\setminus X_k
$$
is a simple homotopy equivalence. However for the triple 
$X_k\subset X_{k-1}\subset X_j$ the conditions on the boundaries
of tubular neighborhoods (2.3) are satisfied. For  such triple
the result was proved in \cite{28, Proposition 2.1} using 
properties of simple homotopy equivalences on triads 
from \cite{8}. \qed
\enddemo
\smallskip
 
 The groups
$LT_*(X,Y,Z)$ and the map
$$
\Theta_*(f,b): [X, G/TOP] \to LT_{n-q-q^{\prime}}(X,Y,Z)
$$
were defined  in  \cite{26} so that the normal map $(f,b)$ is normally bordant
to the  $s$@-triangulation of the triple  $(X,Y,Z)$ if and only if
$\Theta_*(f,b)=0$ (for $n-q-q^{\prime}\geq 5$).

These groups were defined on the spectra level.  First 
we recall necessary facts about application spectra to 
$L$-theory.

 A  spectrum $\Bbb E$ consists of a collection of 
 $CW$--complexes 
 $\{(E_n,*)\}$, $n \in \Bbb Z$, with a collection of cellular 
 maps $\{\epsilon_n:SE_n\to E_{n+1}\}$, 
where $SE_n$ is the suspension of the space $E_n$ \cite{33}. 
The adjoint maps 
$\epsilon_n^{\prime}:E_n\to \Omega E_{n+1}$ (see \cite{33}) 
are defined and the 
$E$ is $\Omega$- spectrum if all 
adjoint maps are homotopy equivalences. 
Let $\Sigma \Bbb E$ be a spectrum  
with  $\{\Sigma \Bbb E\}_n=\Bbb E_{n+1}$ and 
$\{\Sigma\epsilon\}_n=\epsilon_{n+1}$. 
The functor $\Sigma$ has an inverse functor $\Sigma^{-1}$ and  
iterated functors $\Sigma^{k}, k\in \Bbb Z$ 
on the category of spectra are defined. 
For any spectrum $E$ we have an isomorphism 
$$
\pi_n(\Bbb E)=\pi_{n+k}(\Sigma^k \Bbb E)
$$ 
of homotopy groups. 
Recall now that in homotopy theory of spectra there is 
 an equivalence 
between pullback and pushout squares. 
A homotopy commutative square of 
spectra 
$$ 
\matrix 
\Bbb G &\to & \Bbb H \\ 
\downarrow & & \downarrow \\ 
\Bbb E &\to & \Bbb F \\ 
\endmatrix 
\tag 2.7
$$ 
is a  pullback if the 
fibers of horizontal or vertical maps are 
naturally homotopy equivalent  \cite{33}. 
Square (2.7) is a  pushout if the cofibres of
vertical or horizontal maps are naturally homotopy 
equivalent.  

Such natural maps of $L$-groups as transfer and induced map
are realized on the spectra level. 
A homomorphism of oriented groups 
 $f : \pi \to \pi'$ induces 
a cofibration of $\Omega$--spectra (see \cite{11})
$$ 
\CD 
\Bbb L(\pi) @>>> \Bbb L(\pi') @>>> \Bbb L(f) 
\endCD 
\tag 2.8 
$$ 
where $\pi_n(\Bbb L(\pi))=L_n(\pi)$ and 
similarly for the other spectra. 
The homotopy long exact sequence of cofibration (2.8) gives the 
relative 
exact sequence of $L$-groups  
$$ 
\dots\to L_n(\pi)\to L_n(\pi^{\prime})\to L_n(f)\to 
L_{n-1}(\pi)\to\dots. 
$$ 

For a fibration $p : E^{m+n}\to X^n$  over a closed topological 
manifold $X^n$  the transfer map 
$$ 
p^*:L_n(\pi_1(X))\to L_{n+m}(\pi_1(E)) 
$$ 
is defined (see \cite{18}, \cite{19},
 \cite{34}, and \cite{35})
which is realized on the spectra level by a map of 
$\Omega$-spectra 
$$ 
p^{!} : \Bbb L(\pi_1(X)) \to \Sigma^{-m} \Bbb L(\pi_1(E)). 
\tag 2.9
$$ 

For a manifold pair $(X,Y)$ we have the following 
homotopy commutative diagram of spectra 
$$ 
 \matrix 
 \Bbb L(\pi_1(Y)) & \overset{p_1^!}\to{\rightarrow} 
 &\Sigma^{-q}\Bbb 
 L(\pi_1(\partial U)\to \pi_1(U))& 
\overset{\alpha}\to{\rightarrow} & \Sigma^{-q}\Bbb 
 L(\pi_1(X\setminus Y)\to \pi_1(X)) \\ 
 & \ p^!\searrow & \downarrow \delta & & 
 \downarrow\delta_1 \\ 
 & & \Sigma^{1-q}\Bbb L(\pi_1(\partial U))& 
 \overset{\beta}\to{\rightarrow} &\Sigma^{1-q}\Bbb 
 L(\pi_1(X\setminus Y)), 
\endmatrix 
\tag 2.10 
$$ 
 where the left maps are transfer maps and the right 
 horizontal 
 maps 
 are induced by the horizontal maps of the square $F$ (1.6). 
The two right vertical maps in (2.10) obtained 
from extending cofibration sequences (2.8) for vertical maps
of the square $F$ (1.6).
The spectrum $\Bbb LS(F)$ is a homotopy 
cofiber of  the map 
$$ 
 \Sigma^{-1}(\alpha p_1^{!}): \Sigma\Bbb L(\pi_1(Y))\to 
 \Sigma^{-q-1}\Bbb L(\pi_1(X\setminus Y)\to\pi_1(X)) 
\tag 2.11
 $$ 
and the spectrum $\Bbb LP(F)$ is a homotopy cofiber of the 
map 
$$ 
\Sigma^{-1}(\beta p^{!}): \Sigma^{-1}\Bbb L(\pi_1(Y))\to 
\Sigma^{-q}\Bbb L(\pi_1(X\setminus Y)) 
\tag 2.12
$$ 

We have isomorphisms  
(see \cite{1}, \cite{21}, \cite{24}, and \cite{26}) 
$$ 
\pi_n(\Bbb LS(F))\cong LS_n(F), \qquad 
 \pi_n(\Bbb LP(F))\cong LP_n(F).
\tag 2.13 
$$ 

Denote by  $\Cal S_{n+1}(X, Y, \xi)$ the set of 
 concordance classes of $s$@-triangulations of the 
manifold pair $(X, Y, \xi)$ (see \cite{30}). 

For the triple $Z\subset Y \subset X$  of closed topological 
manifolds  consider 
the square of fundamental groups with orientations
for the splitting problem for the manifold pair 
$Z\subset Y$:
$$ 
\Psi= 
\left(\matrix 
\pi_1(\partial V)& \to &\pi_1(Y\setminus Z)\\ 
 \downarrow & & \downarrow\\ 
 \pi_1(Z)&\to &\pi_1(Y)\\ 
 \endmatrix\right). 
\tag 2.14 
$$

Conside the commutative diagram (see \cite{30} and \cite{34})
                  $$ 
                  \matrix 
                  \dots\rightarrow &\Cal S_{n+1}(X, Y, \xi)&\rightarrow & 
                  H_n(X;\bold L_{\bullet})& 
                  \overset{\sigma_1}\to{\rightarrow} &LP_{n-q}(F) &\rightarrow 
                  \cdots \\ 
                  &\downarrow & & 
                  \downarrow & 
                  &\downarrow & \\ 
                  \dots\rightarrow &\Cal S_{n-q+1}(Y)&\rightarrow & 
                  H_{n-q}(Y;\bold L_{\bullet})& 
                  \rightarrow &L_{n-q}(Y) &\rightarrow \cdots\\ 
&\downarrow & & 
                  \downarrow & 
                  &\downarrow= & \\ 
                  \dots\rightarrow &LS_{k}(\Psi)&\rightarrow & 
                  LP_k(\Psi))& 
                  \rightarrow &L_{n-q}(Y) &\rightarrow \cdots\\ 
                  \endmatrix 
                  \tag 2.15 
                  $$ 
 in which $k=n-q-q'$ is a dimension of $Z$, 
 and rows are exact sequences. Observe that the bottom two rows
represent  the diagram (1.7) for the manifold pair $(Y,Z)$. 
Diagram (2.15) is  realized on 
 spectra level (\cite{1} and \cite{26}).

In particular, the composition 
$$ 
 LP_{n-q+1}(F)\to \Cal S_{n+1}(X, Y, \xi)\to 
\Cal S_{n-q+1}(Y)\to LS_{n-q-q'}(\Psi) 
$$ 
of maps from diagram (2.15) is realized 
by a composition $\bold v$ of maps of spectra 
$$ 
\Bbb LP(F)\to \Sigma^{-q}
\Bbb S(X, Y, \xi)\to \Bbb S(Y) \to \Sigma^{q'+1}\Bbb LS(\Psi)
$$
where 
$$ 
\pi_{n}(\Bbb S(X, Y, \xi))= \Cal S_{n}(X, Y, \xi), \
\pi_{n}(\Bbb S(Y))= \Cal S_{n}(Y).  
$$ 

The spectrum $\Bbb LT(X,Y,Z)$ is a homotopy cofiber of the 
map 
$$
\Sigma^{-q'-1}\bold v:\Sigma^{-q'-1}\Bbb LP(F)\to \Bbb LS(\Psi) 
\tag 2.16 
$$
 and by definition  
$LT_n(X,Y,Z)=\pi_n(\Bbb LT(X,Y,Z))$ (see \cite{26}). 
The homotopy long exact sequence of cofibration (2.16) 
gives the  exact sequence 
$$ 
\cdots\to LP_{n-q+1}(F)\to LS_{n-q-q'}(\Psi)\to 
LT_{n-q-q'}(X,Y,Z)\to\cdots 
\tag 2.17 
$$ 
\smallskip 

The triple of manifolds $Z\subset Y\subset X$ is a 
stratified topological manifold (see \cite{3} and \cite{35}) 
which we shall denote
by $\Cal X$. 
Hence the  stratified $L$-groups 
$L^{BQ}(\Cal X)$ of Browder-Quinn are defined. 
These groups are realized on spectra level and we recall 
an inductive definition of these groups from \cite{35, p. 129}
using our notations.  
By Remark 2.3,  the triple $Z\subset Y\subset X$ yields 
a pair of manifolds  with boundaries 
$$
\left(Y\setminus Z,\partial(Y\setminus Z)\right)\subset 
\left(X\setminus Z, \partial(X\setminus Z)\right)
\tag 2.18
$$  
where $\partial(Y\setminus Z)\subset 
\partial(X\setminus Z)$ is a manifold pair which coincides with
natural  decomposition of a boundary 
of a tubular neighborhood of $Z$ in $X$.  
Denote by $F_Z$ the square of fundamental groups 
for splitting problem relative boundary
for the manifold pair (2.18), and by $F_U$ the similar square 
for the closed manifold pair $\partial(Y\setminus Z)\subset 
\partial(X\setminus Z)$.
In fact, the geometric definition of transfer map $p^!$ in 
(2.9) and (2.10) 
for the pair $Z\subset X$
(see \cite{18}, \cite{19}
\cite{30} and \cite{34})  
gives a map 
$$
p^{\#}: L_{n-q-q'}(\pi_1(Z))\to LP_{n-q-1}(F_U)
\tag 2.19
$$ 
which is realized on spectra level (see \cite{35})
by a map of spectra 
$$
{\bold p}^{\#}: \Bbb L(\pi_1(Z))\to \Sigma^{-q'+1} \Bbb LP_{n-q-1}(F_U).
\tag 2.20
$$
Consider the composition of the map ${\bold p}^{\#}$ (2.20)
with the map  of spectra 
$$
{\bold b}:\Sigma^{-q'+1} \Bbb LP_{n-q-1}(F_U)\to 
\Sigma^{-q'+1} \Bbb LP_{n-q-1}(F_Z)
$$ 
which is induced by inclusion of the boundary in (2.18).
We obtain a cofibration  of spectra \cite{35}
$$
{\bold p}^{\#} {\bold b}:\Bbb L(\pi_1(Z))\to 
\Sigma^{-q'+1} \Bbb LP(F_Z)\to \Sigma^{-q-q'+1}\Bbb L^{BQ}(\Cal X) 
\tag 2.21
$$
with a cofiber $\Sigma^{-q-q'+1}\Bbb L^{BQ}(\Cal X)$.
By definition (see \cite{3} and \cite{35})
$$
\pi_n(\Bbb L^{BQ}(\Cal X))= L^{BQ}_n(\Cal X).
$$

For the groups $L^{BQ}_n$  index $n$ 
is equal to dimension of the largest manifold of 
filtration taken  $\bmod \ 4$ (see \cite{3} and \cite{35}). 
 For the case of
surgery obstruction groups $LP_n$ 
Wall and Ranicki 
(see \cite{30} and \cite{34}) used  index $n$ 
which corresponds to dimension of the smallest manifold 
from the pair. Similarly to Wall and Ranicki, 
for the surgery obstruction groups on manifold 
triples  $LT_n$ 
index $n$  is equal to dimension of the bottom manifold of 
the filtration. 
 \smallskip 

\proclaim{Remark 2.6} For the manifold triple $(X,Y,Z)$ 
the homotopy long exact sequence of cofibration (2.21)
gives the following exact sequence of obstruction groups 
$$
\cdots\to L_{n-q-q'}(\pi_1(Z))\to 
LP_{n-q-1}(F_Z)\to  L^{BQ}_{n-1}(\Cal X)\to\cdots 
\tag 2.22
$$
where $\Cal X$ denotes the filtration $Z\subset Y\subset X$.
$\square$
\endproclaim
\bigskip

For the pair $(X,Z)$ we denote  the squares of fundamental 
groups in the splitting problem by $\Phi$. 
The groups $LP_*(\Phi)$ 
fit in the exact sequence (see \cite{30} and \cite{34})
$$
\cdots \to L_n(\pi_1(X\setminus Y)) \to LP_{n-q-q'}(\Phi)\to 
L_{n-q-q'}(\pi_1(Z))\to \cdots
\tag 2.23
$$
which is realized on the spectra level similar to (2.12)
by a cofibration of spectra 
$$
\Bbb LP(\Phi)\to 
\Bbb L(\pi_1(Z))\to \Sigma^{-q-q'+1}\Bbb L(\pi_1(X\setminus Z)).
\tag 2.24
$$
By  \cite{26, Theorem 2} the groups $LT_*$ fit in the
commutative diagram  of exact sequences 
$$ 
\smallmatrix 
\rightarrow & {L}_{n}(C) & \longrightarrow & 
LP_{n-q}(F) & 
 \rightarrow & LS_{k-1}(\Psi)& \rightarrow \cr 
 \ & \nearrow \ \ \ \ \ \ \ \ \searrow & \ & \nearrow \ \ \ \ \ 
                  \ \ 
                  \ 
 \searrow 
  & \ & \nearrow \ \ \ \ \ \ \ \ \searrow & \ \cr 
 \ & \ & LT_k(X,Y,Z)& \ & L_{n-q}(\pi_1(Y)) & \ & \ \cr 
 \ & \searrow \ \ \ \ \ \ \ \ \nearrow & \ & \searrow \ \ \ \ \ 
                  \ \ 
                  \ 
 \nearrow 
 & \ & \searrow \ \ \ \ \ \ \ \ \nearrow & \ \cr 
\rightarrow & LS_{k}(\Psi) & \longrightarrow & 
 LP_{k}(\Psi) & 
\longrightarrow & {L}_{n-1}(C) & \rightarrow, 
 \endsmallmatrix 
\tag 2.25 
$$ 
 where $k=n-q-q'$ and $C=\pi_1(X\setminus Y)$.
 Diagram (2.25) is realized on spectra level and contains
the following exact sequence
$$
\cdots \to L_n(\pi_1(X\setminus Y)) \to LT_{n-q-q'}(X,Y,Z)\to 
LP_{n-q-q'}(\Psi )\to \cdots
\tag 2.26
$$   
Exact sequence (2.26) is realized on the spectra level 
by the cofibration
$$
\Bbb LT(X,Y,Z)\to 
\Bbb LP(\Psi)\to \Sigma^{-q-q'+1}\Bbb L(\pi_1(X\setminus Z)).
\tag 2.27
$$

\proclaim{Proposition 2.7}
There exists the following commutative diagram 

$$
\smallmatrix
&\vdots     &  &\vdots     & & \vdots   & \\
&\downarrow &  &
\downarrow &
 &\downarrow & \\
\dots\rightarrow &L_{n}(\pi_1(X\setminus Y))&\rightarrow &
LT_{k}(X,Y,Z)&
\rightarrow &
LP_{k}(\Psi) &\rightarrow \cdots \\
&\downarrow &  &
\downarrow &
 &\downarrow & \\
\dots\rightarrow & L_{n}(\pi_1(X\setminus Z))&\rightarrow &
LP_{k}(\Phi)&
\rightarrow &L_{k}(\pi_1(Z)) &\rightarrow \cdots\\
&\downarrow &  &\downarrow &&\downarrow & \\
\dots\rightarrow & L_{n}(\pi_1(X\setminus Y)\to \pi_1(X\setminus Z))
&\rightarrow &
LS_{n-q-1}(F_Z)&
\to &
L_{n-q-1}(\pi_1(Y\setminus Z)) &\rightarrow \cdots\\
 &\downarrow &  &\downarrow & &\downarrow & \\
 &\vdots     &  &\vdots     & & \vdots   & \\
\endsmallmatrix
\tag 2.28
 $$
where  $k=n-q-q^{\prime}$. Diagram (2.28) is realized on the spectra level.
All the maps in the square 
$$
\matrix 
LT_*(X,Y,Z) &\to & LP_*(\Psi) \\
\downarrow & & \downarrow \\
LP_*(\Phi) &\to & L_*(\pi_1(Z)) \\
\endmatrix
\tag 2.29
$$
of diagram (2.28) are natural forgetful maps.
The  upper two horizontal rows of the diagram (2.28) 
coincide with exact sequences (2.26) and (2.23).
\endproclaim 

\demo{Proof}
Forgetting the submanifold $Y$ induces the 
natural maps 
$$
LT_*(X,Y,Z)\to LP_*(\Phi) \  \text{and}
\ LP_*(\Psi)\to L_*(\pi_1(Z))
$$ 
which are induced by 
the maps of spectra from (2.25) and (2.27). Similarly to 
(2.25)  the 
forgetful map $LP_*(\Psi)\to L_*(\pi_1(Z))$ is 
realized on the spectra level. 
The forgetful map 
$LT_*(X,Y,Z)\to LP_*(\Phi)$ is realized on the spectra 
level by \cite{28, Theorem 3.5}. This map fits in the following
exact sequence 
$$
\cdots \to LT_{n-q-q'}(X,Y,Z)\to LP_{n-q-q'}(\Phi)\to
LS_{n-q-1}(F_Z)\to\cdots
\tag 2.30
$$
It follows from this that we have 
the following homotopy commutative diagram 
of spectra    
$$
\matrix 
\Bbb LT(X,Y,Z) &\to & \Bbb LP(\Psi) \\
\downarrow & & \downarrow \\
\Bbb LP(\Phi) &\to & \Bbb L(\pi_1(Z)). \\
\endmatrix
\tag 2.31
$$ 

Consider a biinfinite homotopy commutative diagram of spectra 
$$
\smallmatrix 
&\vdots & & \vdots & &\vdots &\\
&\downarrow & & \downarrow & &\downarrow &\\
\to &\Bbb LT(X,Y,Z) &\to & \Bbb LP(\Psi)&\to &\Sigma^{-q^{\prime}-q+1}\Bbb L(\pi_1(X\setminus Y))&\to  \\
&\downarrow & & \downarrow & &\downarrow &\\
\to &\Bbb LP(\Phi) &\to & \Bbb L(\pi_1(Z))&\to &\Sigma^{-q^{\prime}-q+1}\Bbb L(\pi_1(X\setminus Z)) &\to\\
&\downarrow & & \downarrow & &\downarrow &\\
\to &\Sigma^{-q^{\prime}+1}\Bbb LS(F_Z) &\to &\Sigma^{-q^{\prime}+1}\Bbb L(\pi_1(Y\setminus Z)) &\to &\Sigma^{-q^{\prime}-q+1}\Bbb L^{rel} &\to\\
&\downarrow & & \downarrow & &\downarrow &\\
&\vdots & & \vdots & &\vdots & \\
\endsmallmatrix
\tag 2.32
$$
where $\Bbb L^{rel}=\Bbb L(\pi_1(X\setminus Y)\to \pi_1(X\setminus Z))$. 
 This diagram is obtained from 
homotopy commutative diagram (2.31) by consideration of cofibrations given by
all maps of diagram (2.31) (see \cite{21} and \cite{33}). 
Application of $\pi_0$ to (2.32) gives 
commutative diagram (2.28). \qed
\enddemo

\bigskip

We now recall  the following technical result from \cite{21}.

\proclaim{Lemma 2.8} Consider a diagram of spectra 
$$
\matrix
                 & \bullet   &            &  \\
                 &\downarrow &\searrow   &  \\
\bullet\rightarrow      & \bullet   &\rightarrow  \bullet &\\
\searrow                & \downarrow &                     & \\
                         &\bullet &                       &,
\endmatrix
$$
in which the row and  the column are cofibrations.
Then the cofibres of the diagonal  maps are naturally homotopy equivalent.
\endproclaim
\demo{Proof} See \cite{21}. \qed \enddemo
\bigskip

\proclaim{Theorem 2.9} Let $\Cal X$ be a filtration 
$Z\subset Y\subset X$ of topological manifolds,   $n$ 
the   dimension of $X$,
$q$ the codimension of $Y$ in $X$, and $q^{\prime}$ the codimension 
of $Z$ in $Y$. 
We have a homotopy equivalence of 
the spectra 
$$
\Bbb LT(X,Y,Z)\simeq \Sigma^{-q-q'}\Bbb L^{BQ}(\Cal X)
$$ 
and hence  an isomorphism 
$LT_{n-q-q^{\prime}}(X,Y,Z)=L^{BQ}_{n}(\Cal X)$ of surgery obstruction groups 
for $n=0,1,2, 3 \bmod 4$.
\endproclaim 

\demo{Proof}  
It follows by Lemma 2.8 
that the cofibres of the diagonal maps of spectra 
$$
\matrix
\Sigma^{-q^{\prime}}\Bbb L(\pi_1(Y\setminus Z))\to 
\Sigma^{-q^{\prime}+1}\Bbb L(\pi_1(X\setminus Y)),\\
\\
\Bbb LT(X,Y,Z)\to \Bbb L(\pi_1(Z)),\\
\\
\Sigma^{-q^{\prime}+q}\Bbb L(\pi_1(X\setminus Z))
\to \Sigma^{-q^{\prime}+1}\Bbb LS(F_Z)
\endmatrix 
\tag 2.33
$$
 in diagram (2.32)
are naturally homotopy equivalent. 
The map of spectra 
$$
\Sigma \Bbb L(\pi_1(Y\setminus Z))\to \Sigma^{-q} \Bbb L(\pi_1(X\setminus Y))
\tag 2.34
$$
is a realization on spectra level of the transfer 
map for  manifold pair $(X\setminus Z, Y\setminus Z)$ --- this 
 follows from diagram 
(2.28). Hence a cofiber of the first map in (2.33) 
coincides with the spectrum $\Sigma^{q^{\prime}+1}\Bbb LP(F_Z)$. 
Hence the cofiber of the second map in (2.33) coincides with this 
one.  We obtain the following cofibration of spectra
$$
\Bbb LT(X,Y,Z)\to \Bbb L(\pi_1(Z))\to \Sigma^{q^{\prime}+1}\Bbb LP(F_Z).
\tag 2.35
$$
Hence (see \cite{33}) the spectrum $\Bbb LT(X,Y,Z)$
 is defined as homotopical fiber 
of the transfer map
$$
\Bbb L(\pi_1(Z))\to \Sigma^{q^{\prime}+1}\Bbb LP(F_Z).
\tag 2.36
$$
However, by (2.21) a homotopical fiber of this map 
is a spectrum
$\Sigma^{-q-q^{\prime}}\Bbb L^{BQ}(\Cal X))$
where $\Cal X$ is the filtration $Z\subset Y\subset X$.
Therefore the assertion of the theorem
follows.
\qed \enddemo
\smallskip

\proclaim{Corollary 2.10} Under hypothesis of Theorem 2.9 
we have the following three 
braids of exact sequences
$$ 
\matrix 
\rightarrow & {L}_{n-q}(D) & \longrightarrow & 
L_{n-1}(C) & 
 \rightarrow & LT_{k-1}& \rightarrow \cr 
 \ & \nearrow \ \ \ \ \ \ \ \ \searrow & \ & \nearrow \ \ \ \ \ 
                  \ \ 
                  \ 
 \searrow 
  & \ & \nearrow \ \ \ \ \ \ \ \ \searrow & \ \cr 
 \ & \ & LP_k(\Psi)& \ & LP_{m}(F_Z) & \ & \ \cr 
 \ & \searrow \ \ \ \ \ \ \ \ \nearrow & \ & \searrow \ \ \ \ \ 
                  \ \ 
                  \ 
 \nearrow 
 & \ & \searrow \ \ \ \ \ \ \ \ \nearrow & \ \cr 
\rightarrow & LT_{k} & \longrightarrow & 
 L_{k}(\pi_1(Z)) & 
\longrightarrow & {L}_{m}(D) & \rightarrow, 
 \endmatrix 
\tag 2.37 
$$ 
$$ 
\matrix 
\rightarrow & LT_{k} & \longrightarrow & 
L_{k}(\pi_1(Z)) & 
 \rightarrow & L_{n-1}(E)& \rightarrow \cr 
 \ & \nearrow \ \ \ \ \ \ \ \ \searrow & \ & \nearrow \ \ \ \ \ 
                  \ \ 
                  \ 
 \searrow 
  & \ & \nearrow \ \ \ \ \ \ \ \ \searrow & \ \cr 
 \ & \ & LP_k(\Phi)& \ & LP_{m}(F_Z) & \ & \ \cr 
 \ & \searrow \ \ \ \ \ \ \ \ \nearrow & \ & \searrow \ \ \ \ \ 
                  \ \ 
                  \ 
 \nearrow 
 & \ & \searrow \ \ \ \ \ \ \ \ \nearrow & \ \cr 
\rightarrow & L_{n}(E) & \longrightarrow & 
 LS_{m}(F_Z) & 
\longrightarrow & {LT}_{k-1} & \rightarrow, 
 \endmatrix 
\tag 2.38 
$$ 
and
$$ 
\matrix 
\rightarrow & L_{n}(E) & \longrightarrow & 
LS_{m}(F_Z) & 
 \rightarrow & L_{m}(D)& \rightarrow \cr 
 \ & \nearrow \ \ \ \ \ \ \ \ \searrow & \ & \nearrow \ \ \ \ \ 
                  \ \ 
                  \ 
 \searrow 
  & \ & \nearrow \ \ \ \ \ \ \ \ \searrow & \ \cr 
 \ & \ & L_n(C\to E)& \ & LP_{m}(F_Z) & \ & \ \cr 
 \ & \searrow \ \ \ \ \ \ \ \ \nearrow & \ & \searrow \ \ \ \ \ 
                  \ \ 
                  \ 
 \nearrow 
 & \ & \searrow \ \ \ \ \ \ \ \ \nearrow & \ \cr 
\rightarrow & L_{n-q}(D) & \longrightarrow & 
 L_{n-1}(C) & 
\longrightarrow & L_{n-1}(E) & \rightarrow, 
 \endmatrix 
\tag 2.39 
$$ 
 where $k=n-q-q'$, $m=n-q-1$, $D=\pi_1(Y\setminus Z)$, 
 $E=\pi_1(X\setminus Z)$ and $C=\pi_1(X\setminus Y)$.
 Diagrams (2.37), (2.38), and  (2.39) are realized on the spectra level.
\endproclaim

\demo{Proof }
From biinfinite homotopy commutative diagram  (2.32) and cofibration 
(2.35) we obtain the folllowing  homotopy  commutative diagram of
spectra
$$
\CD
\Bbb LT(X,Y,Z) @>>> \Bbb LP(\Psi)
@>>>
 \Sigma^{-q^{\prime}-q+1}\Bbb L(\pi_1(X\setminus Y))
\\
 @V=VV  @VVV  @VVV \\
\Bbb LT(X,Y,Z) @>>> \Bbb L(\pi_1(Z))
@>>>
 \Sigma^{q^{\prime}+1}\Bbb LP(F_Z)\\
\endCD
\tag 2.40
$$
in which horizontal rows  are cofibrations, and right vertical map 
is induced by two left vertical maps (see \cite{33}). Hence fibers
of the two right horizontal maps in (2.40) are naturally 
homotopy equivalent to the spectrum $\Bbb LT(X,Y,Z)$. Hence the 
right square in
(2.40) is a pullback and fibers of vertical map of this square 
are also naturally homotopy equivalent. Homotopy long exact 
sequences of this square give commutative diagram (2.37). In 
a similar way the  
commutative diagrams  (2.38) and (2.39) follow from 
the other two cofibrations from (2.33) and homotopy 
commutative diagram (2.32).
\qed
\enddemo
\smallskip

\proclaim{Remark 2.11} Diagram (2.39) is, in fact, diagram  (1.8) constructed
for the pair of manifolds with boundaries
$(X\setminus Y)\subset (X\setminus Z)$.
\qed 
\endproclaim

\proclaim{Remark 2.12} We can consider a manifold pair 
$Y^{n-q}\subset X^n$ as the stratified manifold $\Cal X$ for which 
Browder-Quinn 
groups $L^{BQ}(\Cal X)$ are defined.
It  follows from cofibration (2.12), that 
the definition of $LP_*$-groups by Wall and results 
of Ranicki (see \cite{30} and \cite{34}) 
give an isomorphism $LP_{n-q}(F)\cong L^{BQ}(\Cal X)$. This isomorphism 
is realized on the spectra level. 
\qed 
\endproclaim 

\bigskip

\subhead 3. Surgery on a manifold with filtration 
\endsubhead
\bigskip

In  this section we introduce surgery obstruction groups 
for the filtration $\Cal X$ (1.11) and describe theier main 
properties. At first we give the motivation of our definition 
and then we prove Theorem 3.1 and describe relations of introduced groups
to $L^{BQ}$-groups of Browder and Quinn.  We shall use 
the notations of 
previous sections. 

For a manifold pair $(X^n, Y^{n-q})$ of codimension $q$ realization of 
the diagram (1.8) on spectra level provides the 
following homotopy commutative diagram 
of spectra  
$$
\matrix
               &              & \Sigma^q\Bbb L(\pi_1(Y))               &                 &                                 \\
                &             &\downarrow                              &\searrow            &                                  \\
\Bbb L(\pi_1(X))& \rightarrow& \Bbb L(\pi_1(X\setminus Y)\to \pi_1(X)) &\rightarrow     &\Sigma \Bbb L(\pi_1(X\setminus Y))\\
                & \searrow& \downarrow                            &                &                                \\
                &         &\Sigma^{q+1}\Bbb LS(F),                      &                &                               
\endmatrix
\tag 3.1
$$
in which the vertical column and horizontal row are cofibrations.
The cofibres of diagonal maps are naturally homotopy equivalent to 
the spectrum $\Sigma^{q+1} \Bbb LP(F)$ as follows from (1.8) and Lemma 2.8.  

Consider now a manifold triple  
$Z^{n-q-q^{\prime}}\subset Y^{n-q} \subset X^n$ where $q$  is 
the codimension of $Y$ in $X$ and $q^{\prime}$ is the codimension $Z$ in $Y$. 
Realization of 
the diagram (2.25) on spectra level provides the 
following homotopy commutative diagram 
of spectra  
$$
\matrix
               &              & \Sigma^{q^{\prime}}\Bbb LP(\Psi)  &                 &                                 \\
                &             &\downarrow                              &\searrow            &                                  \\
\Bbb LP(F)& \rightarrow& \Bbb L(\pi_1(Y)) &\rightarrow     &\Sigma^{-q+1} \Bbb L(\pi_1(X\setminus Y))\\
                & \searrow& \downarrow                            &                &                                \\
                &         &\Sigma^{q^{\prime}+1}\Bbb LS(\Psi),                      &                &                               
\endmatrix
\tag 3.2
$$
in which the vertical column and horizontal row are cofibrations.
Recall that $F$ is a square of fundamental groups  for 
splitting problem for the pair $(X,Y)$, and $\Psi$ is a similar
square for the pair $(Y,Z)$. 
The cofibres of the diagonal  maps are naturally homotopy equivalent to 
the spectrum $\Sigma^{q^{\prime}+1} \Bbb LT(X,Y,Z)$ 
as follows from (2.25) and Lemma 2.8.  

Now consider filtration $\Cal X$ (1.11) for which the restricted 
filtrations $\Cal X_j$ for $j=0, 1, ..., k$ are defined. 

For a pair of submanifolds $X_{j}\subset X_{j-1}$ of filtration (1.11) we 
denote the square of fundamental groups for splitting problem by 
$F_j$ where $1\leq j\leq k$. 
We also introduce special notations for the following 
filtrations. Let $\Cal Y$ be a subfiltration 
$$
X_k\subset X_{k-1}\subset \cdots \subset X_2\subset X_1
\tag 3.3
$$
of $\Cal X$
and $\Cal Y_{j-1}$ be a restricted subfiltration 
$$
X_j\subset X_{j-1}\subset \cdots \subset X_2\subset X_1
\tag 3.4
$$
of $\Cal Y$
where $1\leq j\leq k$.
 We have $\Cal X_0 = (X_0)=(X)$, $\Cal X_1=(X_1\subset X_0)$, 
and $\Cal X_2=(X_2\subset X_1\subset X_0)$. 
Denote 
$$
\Bbb LM^0(\Cal X)=\Bbb LM^0(X_0)=\Bbb L(\pi_1(X_0)),
$$   
$$
\Bbb LM^1(\Cal X)=\Bbb LM^1(\Cal X_1)=\Bbb LM^1(X_1\subset X_0)=
\Bbb LP(F_1),
$$
and 
$$
\Bbb LM^2(\Cal X)=\Bbb LM^2(\Cal X_2)=\Bbb LT(X_0,X_1, X_2).
$$
For  spectra $\Bbb LM^i$ defined above 
 with $0\leq i\leq 2$ and  $j\geq i$
we have, by definition, that  $\Bbb LM^i(\Cal X_j)=\Bbb LM^i(\Cal X)$.
Diagram (3.2) 
in our notations has the following form 
$$
\matrix
               &              & \Sigma^{q_2}\Bbb LM^1(\Cal Y)  &                 &                                 \\
                &             &\downarrow                              &\searrow            &                                  \\
\Bbb LM^1(\Cal X)& \rightarrow& \Bbb LM^0(\Cal Y) &\rightarrow     &\Sigma^{-q_1+1} \Bbb L(\pi_1(X_0\setminus X_1))\\
                & \searrow& \downarrow                            &                &                                \\
                &         &\Sigma^{q_2+1}\Bbb LS(F_2)                      &                &                               
\endmatrix
\tag 3.5
$$
with the  cofibres of diagonal maps which are naturally 
homotopy equivalent to 
$$
\Sigma^{q_2+1} \Bbb LM^2(\Cal X)=\Sigma^{q_2+1} \Bbb LM^2(\Cal X_2)
$$

 The right diagonal map 
from diagram (3.5) gives a cofibration of spectra 
$$
\Bbb LM^2(\Cal X)\to \Bbb LM^1(\Cal Y) \to 
\Sigma^{-q_1-q_2+1}\Bbb L(\pi_1(X_0\setminus X_1))
\tag 3.6
$$
where $-q_1-q_2=l_2-n$.
The left diagonal map from (3.5) 
gives a cofibration 
$$
\Sigma^{q_2}\Bbb LM^2(\Cal X)\to \Bbb LM^1(\Cal X)\to 
\Sigma^{q_2+1}\Bbb LS(F_2).
\tag 3.7
$$
For the filtration $\Cal Y$ cofibration (3.7) gives a cofibration
$$
\Sigma^{q_3}\Bbb LM^2(\Cal Y)\to \Bbb LM^1(\Cal Y)\to 
\Sigma^{q_3+1}\Bbb LS(F_3).
\tag 3.8
$$
We can combine cofibrations (3.6) and (3.8) to obtain the following  
homotopy commutative diagram 
$$
\matrix
               &              & \Sigma^{q_3}\Bbb LM^2(\Cal Y)  &                 &                                 \\
                &             &\downarrow                              &\searrow            &                                  \\
\Bbb LM^2(\Cal X)& \rightarrow& \Bbb LM^1(\Cal Y) &\rightarrow     &\Sigma^{l_2-n+1} \Bbb L(\pi_1(X_0\setminus X_1))\\
                & \searrow& \downarrow                            &                &                                \\
                &         &\Sigma^{q_3+1}\Bbb LS(F_3)                      &                &                               
\endmatrix
\tag 3.9
$$
in which cofibers of diagonal maps are naturally homotopy equivalent.
We shall denote homotopy cofiber of diagonal map in diagram (3.9) 
by 
$$
\Sigma^{q_3+1}\Bbb LM^3(\Cal X_3) =\Sigma^{q_3+1}\Bbb LM^3(\Cal X).
$$ 
It follows from this definition that 
$\Bbb LM^3(\Cal X_j)=\Bbb LM_3(\Cal X)$ for $3\leq j\leq k$. 
We can continue these constructions to give inductive definition of 
the spectra 
$$
\Bbb LM^i(\Cal X) =\Bbb LM^i(\Cal X_i)
$$
for $4\leq i\leq k$.

Let a spectrum $\Bbb LM^j(\Cal X)=\Bbb LM^j(\Cal X_j)$ be already defined
for $k\geq j\geq 2$ in such a  way that the 
spectrum $\Sigma^{q_j+1}\Bbb LM^j(\Cal X),\ (j\geq 2) $ 
is a natural homotopy cofiber    of diagonal maps in a diagram 
$$
\matrix
               &              & \Sigma^{q_j}\Bbb LM^{j-1}(\Cal Y)  &                 &                                 \\
                &             &\downarrow                              &\searrow            &                                  \\
\Bbb LM^{j-1}(\Cal X)& \rightarrow& \Bbb LM^{j-2}(\Cal Y) &\rightarrow     &\Sigma^{l_{j-1}-n+1} \Bbb L(\pi_1(X_0\setminus X_1))\\
                & \searrow& \downarrow                            &                &                                \\
                &         &\Sigma^{q_j+1}\Bbb LS(F_j).                      &                &                               
\endmatrix
\tag 3.10
$$
 The right diagonal map
from diagram (3.10) gives a cofibration of the spectra 
$$
\Bbb LM^{j}(\Cal X)\to \Bbb LM^{j-1}(\Cal Y) \to 
\Sigma^{l_{j-1}-n-q_j+1}\Bbb L(\pi_1(X_0\setminus X_1))
\tag 3.11
$$
where $l_{j-1}-n-q_j+1=l_{j}-n+1$.
The left diagonal map in (3.10) 
gives the cofibration 
$$
\Sigma^{q_j}\Bbb LM^j(\Cal X)\to \Bbb LM^{j-1}(\Cal X)\to 
\Sigma^{q_j+1}\Bbb LS(F_j).
\tag 3.12
$$
For the filtrations $\Cal Y$ and $\Cal Y_j$ 
 cofibration (3.12) gives the cofibration
$$
\Sigma^{q_{j+1}}\Bbb LM^j(\Cal Y)\to \Bbb LM^{j-1}(\Cal Y)\to 
\Sigma^{q_{j+1}+1}\Bbb LS(F_{j+1}).
\tag 3.13
$$
We can combine cofibrations (3.11) and (3.13) to obtain the following 
homotopy commutative diagram 
$$
\matrix
               &              & \Sigma^{q_{j+1}}\Bbb LM^j(\Cal Y)  &                 &                                 \\
                &             &\downarrow                              &\searrow            &                                  \\
\Bbb LM^j(\Cal X)& \rightarrow& \Bbb LM^{j-1}(\Cal Y) &\rightarrow     &\Sigma^{l_j-n+1} \Bbb L(\pi_1(X_0\setminus X_1))\\
                & \searrow& \downarrow                            &                &                                \\
                &         &\Sigma^{q_{j+1}+1}\Bbb LS(F_{j+1}).                      &                &                               
\endmatrix
\tag 3.14
$$
in which cofibers of diagonal maps are naturally homotopy equivalent.
We shall denote homotopy cofiber of the diagonal map in diagram (3.14) 
by 
$$
\Sigma^{q_{j+1}+1}\Bbb LM^{j+1}(\Cal X) =\Sigma^{q_{j+1}+1}\Bbb LM^j(\Cal X_j).
$$

Thus for $0\leq j\leq k$ the spectra $\Bbb LM^j(\Cal X)$ are defined.
It follows from the definition that 
$$
\Bbb LM^j(\Cal X)= \Bbb LM^j(\Cal X_i), \ \text{for} \ k\geq i\geq j.
$$ 
We define groups $LM^i_j(\Cal X)$ as homotopy groups
$\pi_j(\Bbb LM^j(\Cal X)$. It follows from definition, 
that $j$ is defined by $\bmod \ 4$.  
\smallskip

\proclaim{Proposition 3.1} Let $\Cal X$ be  filtration (1.11). 
For $0\leq i\leq k-2$ the groups $LM$ fit in the following braid 
of exact sequences
$$ 
\matrix 
\rightarrow & LS_{l_j}(F_j) & \longrightarrow & 
LM^{j-1}_{l_j}(\Cal Y) & 
 \rightarrow & L_{n-1}(\pi_1(X_0\setminus X_1))& \rightarrow \cr 
 \ & \nearrow \ \ \ \ \ \ \ \ \searrow & \ & \nearrow \ \ \ \ \ 
                  \ \ 
                  \ 
 \searrow 
  & \ & \nearrow \ \ \ \ \ \ \ \ \searrow & \ \cr 
 \ & \ & LM^j_{l_j}(\Cal X)& \ & LM_{l_{j-1}}^{j-2}(\Cal Y) & \ & \ \cr 
 \ & \searrow \ \ \ \ \ \ \ \ \nearrow & \ & \searrow \ \ \ \ \ 
                  \ \ 
                  \ 
 \nearrow 
 & \ & \searrow \ \ \ \ \ \ \ \ \nearrow & \ \cr 
\rightarrow & L_{n}(\pi_1(X_0\setminus X_1)) & \rightarrow & 
 LM^{j-1}_{l_{j-1}}(\Cal X) & 
\rightarrow & LS_{l_j-1}(F_j) & \rightarrow, 
 \endmatrix 
\tag 3.15 
$$ 
 where $l_j$ is the dimension of the bottom manifold of filtration. 
 This diagram is realized on the spectra level.
 \endproclaim
\demo{Proof} From the definition of $LM$-groups we get 
a homotopy commutative square of the spectra 
$$
\matrix
\Bbb LM^{j-2}(\Cal Y) &\to     
&\Sigma^{l_{j-1}-n+1} \Bbb L(\pi_1(X_0\setminus X_1))\\
\downarrow   &                &  \downarrow \\
\Sigma^{q_j+1}\Bbb LS(F_j)&\to & \Sigma^{q_j+1}\Bbb LM^j(\Cal X).\\
\endmatrix
\tag 3.16
$$
The fibres of parallel maps in (3.16) are naturally homotopy equivalent --- 
this follows from diagram (3.10). Hence square (3.16) is a pullback and 
consideration of homotopy long exact sequences of the maps from 
this square completes the proof of the theorem.   
\qed
\enddemo
\smallskip

\proclaim{Corollary 3.2} For $ 2\leq j\leq k$ 
the spectrum $\Bbb LM^j(\Cal X)$ fits
in the following pullback square of spectra
$$
\matrix
\Sigma^{q_j}\Bbb LM^{j}(\Cal X) &\to     
&\Sigma^{q_{j}} \Bbb LM^{j-1}(\Cal Y)\\
\downarrow   &                &  \downarrow \\
\Bbb LM^{j-1}(\Cal X)&\to & \Bbb LM^{j-2}(\Cal Y).\\
\endmatrix
\tag 3.17
$$ 
\qed
\endproclaim

We can now  define the spectra for the structure sets of filtration 
$\Cal X$. 
In accordance with Ranicki \cite{30}  we define a spectrum 
$\Bbb S^0(\Cal X)=\Bbb S(X)$ for a manifold $X^n$
as a homotopical fiber of the  map (1.4). 

For a filtration 
$\Cal X$ which is given by a manifold pair 
$(X^n, Y^{n-q})$ the map 
$$
H_n(X;\bold L_{\bullet})\to LP_{n-q}(F)
$$ 
from (1.15) is 
realized on the 
spectra level by a map of spectra (see \cite{29}, \cite{30}, \cite{1},
and \cite{26}) 
$$
X_+\land \bold L_{\bullet}\to \Sigma^{q}\Bbb LP(F).
\tag 3.18
$$
We denote the cofiber of the map 
in (3.18) by $\Bbb S^1(\Cal X)=\Bbb S(X,Y, \xi)$ 
with homotopy groups 
$$
\pi_n(\Bbb S^1(\Cal X))=S_{n}(X, Y,\xi)
$$
fitting in the exact sequence (1.15). 

For a filtration 
$\Cal X=(Z\subset Y\subset X)$ we define the spectrum 
$\Bbb S^2(\Cal X)= \Bbb S(X,Y,Z)$ (see \cite{26}) as a
homotopical cofiber of the map
$$
X_+\land \bold L_{\bullet}\to \Sigma^{q+q^{\prime}}\Bbb LT(X,Y,Z).
\tag 3.19
$$   

A $t$-triangulation of filtration (1.11) $\Cal X$ gives $t$-triangulations
of restricted filtrations $\Cal X_k$,
$\Cal Y$, and $\Cal Y_{k-1}$. Thus we obtain the following
 commutative diagram  
$$
\matrix
\Cal T(\Cal X)&\to &\Cal T(\Cal Y)\\
\downarrow& & \downarrow\\
\Cal T(\Cal X_k)&\to & \Cal T(\Cal Y_k)\\
\endmatrix
\tag 3.20
$$
which is realized on the spectra level (see \cite{29}, \cite{30}, and \cite{26}).
By Proposition 2.4 (see \cite{30}), the diagram (3.20) 
on the spectra level has the following form
$$
\Cal F=\left(\matrix
(X_0)_+\land\bold L_{\bullet}&\to &
\Sigma^{q_1}[(X_1)_+\land\bold L_{\bullet}]\\
\downarrow& & \downarrow\\
(X_0)_+\land\bold L_{\bullet}&\to &
\Sigma^{q_1}[(X_1)_+\land\bold L_{\bullet}]\\
\endmatrix\right).
\tag 3.21
$$

It follows from the definition of spectra $\Bbb LM^j(\Cal X)$ and from 
(1.4), (3.18) and (3.19) that for $k\geq j\geq 0$ we have 
the maps 
$$
(X_0)_+\land \bold L_{\bullet}\to \Sigma^{n-l_j}\Bbb LM^j(\Cal X)
\tag 3.22
$$
with cofiber which we shall denote by $\Bbb S^j(\Cal X)$.
Thus  
$\Bbb S^0(\Cal X)=\Bbb S(X_0)$.
Using the maps from (3.22) 
we obtain a map $\Lambda_2$ of squares 
$$
\Lambda_2 :\Cal F\longrightarrow \Cal G_2
\tag 3.23
$$
where 
$$
\Cal G_2=
\left(\matrix
\Sigma^{n-l_2}\Bbb LM^{2}(\Cal X) &\to     
&\Sigma^{n-l_2} \Bbb LM^{1}(\Cal Y)\\
\downarrow   &                &  \downarrow \\
\Sigma^{n-l_1}\Bbb LM^{1}(\Cal X)&\to & \Sigma^{n-l_1} \Bbb LM^{0}(\Cal Y).\\
\endmatrix\right)
\tag 3.24
$$ 
which gives a homotopy commutative diagram of spectra in  form
of a cube. Observe  here that  square in (3.24) 
follows from (3.17). 

The cofibres of four maps which constitute the map  $\Lambda_2$
give a pullback square 
$$
\matrix
\Bbb S^2(\Cal X)&\to &\Sigma^{q_1}\Bbb S^1 (\Cal Y)\\
\downarrow& & \downarrow\\
\Bbb S^1(\Cal X)&\to &\Sigma^{q_1}\Bbb S^0(\Cal Y).\\
\endmatrix
\tag 3.25
$$
since squares (3.17) and (3.21) are pullback.

Let   $\Cal G_i, \ k\geq i\geq 2$ be a homotopy commutative
square of spectra 
$$
\Cal G_i=
\left(\matrix
\Sigma^{n-l_i}\Bbb LM^{i}(\Cal X) &\to     
&\Sigma^{n-l_{i}} \Bbb LM^{i-1}(\Cal Y)\\
\downarrow   &                &  \downarrow \\
\Sigma^{n-l_{i-1}}\Bbb LM^{i-1}(\Cal X)&\to & 
\Sigma^{n-l_{i-1}} \Bbb LM^{i-2}(\Cal Y)\\
\endmatrix\right)
\tag 3.26
$$
which 
follows from (3.17). 

\proclaim{Proposition 3.3} 
For $k\geq i\geq 2$
there exist  maps 
$$
\Lambda_i:\Cal F\longrightarrow \Cal G_i
\tag 3.27
$$
of squares which are given by four maps in such a way that 
the resulting  diagram in form of a cube is homotopy commutative. 
\endproclaim

\demo{Proof} Using induction on  $i$ it suffices 
to define the left upper map in $\Lambda_i$  
when the other three maps
are already defined.  
This is possible since homotopy commutative square
(3.26) is a pullback. \qed
\enddemo
\smallskip

\proclaim{Definition} Let $\Cal X$ be a filtration (1.11). 
For $k\geq j\geq 0$ we shall denote by
$\Bbb S^j(\Cal X)$ a homotopical cofiber of the map
$$
(X_0)_+\land \bold L_{\bullet}\to \Sigma^{n-l_j} \Bbb LM^j(\Cal X)
\tag 3.28
$$
which is given by the map $\Lambda_i$ in (3.27). 
We shall denote the homotopy groups  
$\pi_n(\Bbb S^j(\Cal X))$  by
$\Cal S^j_n(\Cal X)$.\endproclaim 
\smallskip

The structure sets $\Cal S^j_n(\Cal X)$ are  natural generalizations 
of the structure sets of homotopy triangulations $\Cal S_n(X)$ 
of a manifold $X$ and homotopy triangulations of a manifold 
pair $\Cal S_n(X,Y,\xi)$. Now we shall describe the main properties
of the introduced sets.
\smallskip

\proclaim{Remark 3.4} Let $\Cal X$ be the filtration (1.11). 
It  follows from definition for $k\geq j\geq 0$  
that we have the exact sequence
$$  
\cdots\to \Cal S^j_{n+1}(\Cal X)\to H_n(X;\bold L_{\bullet})\to
LM^j_{l_j}(\Cal X)\to\cdots .
\tag 3.29
$$
For $j=0$ the exact sequence (3.29)  coincides with (1.2) with 
$X=X_0$, for $j=1$ 
it coincides with (1.15) for the pair $X_1\subset X_0$, 
and for $j=2$  it coincides 
with homotopy long exact sequence of cofibration (3.19)
for the triple $X_2\subset X_1\subset X_0$. \qed
\endproclaim 
\smallskip

\proclaim{Proposition 3.5} For $k\geq i\geq 2$ 
there  exist the following homotopy commutative 
pullback squares of spectra
$$
\matrix
\Bbb S^i(\Cal X)&\to &\Sigma^{q_1}\Bbb S^{i-1} (\Cal Y)\\
\downarrow& & \downarrow\\
\Bbb S^{i-1}(\Cal X)&\to &\Sigma^{q_1}\Bbb S^{i-2}(\Cal Y).\\
\endmatrix
\tag 3.30
$$
\endproclaim
\demo{Proof} The square (3.30) obtained as square of homotopical cofibres 
of the maps constitute the map $\Lambda_i$. 
Squares (3.21) and (3.27) are pullback. Hence square (3.30) is 
a pullback. \qed
\enddemo
\smallskip

\proclaim{Corollary 3.6} For $k\geq i\geq 2$ there exist the following
braids of exact sequences
$$ 
\matrix 
\rightarrow & \Cal S_{n}(X_0\setminus X_1) & 
\longrightarrow & 
\Cal S_n^{i-1}(\Cal X) & 
         \rightarrow & LS_{l_i-1}(F_i)& \rightarrow \cr 
    \ & \nearrow \ \ \ \ \ \ \ \ \searrow & \ & \nearrow \ \ \ \ \ 
    \ \    \  \searrow 
 & \ & \nearrow \ \ \ \ \ \ \ \ \searrow & \ \cr 
\ & \ & \Cal S_{n}^i(\Cal X)& \ & \Cal S_{n-q_1}^{i-2}(\Cal Y) & \ & \ \cr 
                  \ & \searrow \ \ \ \ \ \ \ \ \nearrow & \ & \searrow \ \ \ \ \ 
                  \ \ 
                  \ 
                  \nearrow 
                  & \ & \searrow \ \ \ \ \ \ \ \ \nearrow & \ \cr 
       \rightarrow & LS_{l_i}(F_i) & \longrightarrow & 
        \Cal S^{i-1}_{n-q_1}(\Cal Y) & 
 \longrightarrow & \Cal S_{n-1}(X_0\setminus X_1) & \rightarrow 
\endmatrix 
\tag 3.31
$$ 
where $\Cal S_n(X_0\setminus X_1)$ is the structure set fitting 
in  algebraic surgery exact sequence (1.2) for $X_0\setminus X_1$.                  
\endproclaim 
\demo{Proof} The homotopy long exact sequences of the maps from the pullback 
square (3.30) give the diagram (3.31). \qed
\enddemo
\smallskip

Diagram (3.31) for the case of a manifold 
triple  was obtained in 
 \cite{26}. In fact, this diagram is a natural generalization of 
the diagram \cite{30, Proposition 7.2.6 ii)} which was given there for
manifold pairs.

\smallskip

\proclaim{Proposition 3.7} For $k\geq i\geq 1$ 
there  exist the following homotopy commutative 
pullback squares of spectra
$$
\matrix
\Sigma^{-1}\Bbb S^{i-1}(\Cal X)&\to &(X_0)_+\land \bold L_{\bullet}\\
\downarrow& & \downarrow\\
\Sigma^{n-l_i} \Bbb LS(F_i)&\to &\Sigma^{n-l_i}\Bbb LM^i(\Cal X).\\
\endmatrix
\tag 3.32
$$
\endproclaim
\demo{Proof} We have the following homotopy commutative 
diagram 
$$
\matrix
\Sigma^{-1}\Bbb S^{i-1}(\Cal X)&\to &(X_0)_+\land \bold L_{\bullet}&\to &
\Sigma^{n-l_{i-1}}\Bbb LM^{i-1}(\Cal X)\\
\downarrow& & \downarrow & & \downarrow=\\
\Sigma^{n-l_i} \Bbb LS(F_i)&\to &\Sigma^{n-l_i}\Bbb LM^i(\Cal X)
&\to& \Sigma^{n-l_{i-1}}\Bbb LM^{i-1}(\Cal X)\\
\endmatrix
\tag 3.33
$$
where the right square follows from Proposition 3.3 and the left map 
is obtained as a natural map of fibres by \cite{33}. We have 
such fibres  by definition of the spectra $\Bbb S^i(\Cal X)$ 
and by Corollary 3.6. Now the cofibres of horizontal maps
in the left square of 3.33 are naturally homotopy equivalent and 
this square is a  pushout, and hence it is a  pullback.
\qed
\enddemo 
 \smallskip

\proclaim{Corollary 3.8} For $k\geq i\geq 1$ there exist the following
braid of exact sequences
$$ 
\matrix 
\rightarrow & \Cal S_{n+1}^i(\Cal X) & 
\longrightarrow & 
H_n(X_0;\bold L_{\bullet}) & 
         \rightarrow & LM^{i-1}_{l_{i-1}}(\Cal X)& \rightarrow \cr 
    \ & \nearrow \ \ \ \ \ \ \ \ \searrow & \ & \nearrow \ \ \ \ \ 
    \       \searrow\Theta_i 
 & \ & \nearrow \ \ \ \ \ \ \ \ \searrow & \ \cr 
\ & \ & \Cal S_{n+1}^{i-1}(\Cal X)& \ & LM^{i}_{l_i}(\Cal X) & \ & \ \cr 
                  \ & \searrow \ \ \ \ \ \ \ \ \nearrow & \ & \searrow \ \ \ \ \ 
                  \ \ 
                  \ 
                  \nearrow 
                  & \ & \searrow \ \ \ \ \ \ \ \ \nearrow & \ \cr 
       \rightarrow & LM_{l_{i-1}+1}^{i-1}(\Cal X) & \longrightarrow & 
        LS_{l_i}(F_i) & 
 \longrightarrow & \Cal S_{n}^i(\Cal X) & \rightarrow  .
\endmatrix 
\tag 3.34
$$ 
\endproclaim 
\demo{Proof} The homotopy long exact sequences of the maps from pullback 
square (3.32) give commutative diagram of exact sequences (3.31). \qed
\enddemo
\smallskip

Diagram (3.34) for the case of a manifold pair 
 was obtained in 
 \cite{30, Proposition 7.2.6 iv)}. For the case of a manifold triple 
 this diagram was obtained in \cite{26, Theorem 4}.
 
\bigskip

\proclaim{Theorem 3.9}  Let $\Cal X$ be the filtration (1.11) where
 dimension 
of the submanifold $X_k$ is equal to $l_k\geq 5$,  and let 
$$
x=(f,b)\in [X_0, G/TOP]= H_n(X; \bold L_{\bullet})
$$ 
be a $t$-triangulation of $\Cal X$  with the given map
$f: M\to X=X_0$. Then the map $(f,b)$ is normally bordant to 
an $s$@-triangulation of the filtration $\Cal X$  
if and only if 
$\Theta_k(x)=0$. We can identify the set $\Cal S_{n+1}^i(\Cal X)$ 
with the set of concordance classes of $s$-triangulations 
of $\Cal X_i$.
\endproclaim 
\demo{Proof} We use unduction on the number of submanifolds
$k$. For $k=1,2$ the result was obtained in \cite{30} and \cite{26},
respectively.  

Let $x=[(f,b)]\in H_n(X_0;\bold L_{\bullet})$ 
be an $s$-triangulation of filtration $\Cal X_k$. It follows 
from definition that it is an $s$-triangulation of the
subfiltration $\Cal X_{k-1}$ for which the restriction on
$X_{k-1}$  is already split along the submanifold $X_k\subset X_{k-1}$.
Hence, by inductive hypothesis, $\Theta_{i-1}(x)=0$ and it follows 
from (3.34) 
that $x$ represents an element  
$y\in \Cal S_{n+1}^{k-1}(\Cal X_{k-1})$. 
It follows from diagram (3.31) 
that the map 
$$
\sigma:\Cal S^{k-1}_{n+1}(\Cal X)=\Cal S^{k-1}(\Cal X_{k-1})\to LS_{l_k}(F_k)
$$ 
in diagram (3.34) 
is given by the composition 
$$
\Cal S^{k-1}_{n+1}(\Cal X_{k-1})
\to \Cal S_{l_{k-1}+1}(X_{k-1})\to LS_{l_k}(F_k).
\tag 3.35
$$
The last map in (3.35) is the map from diagram (1.7) for the the pair 
$(X_{k-1}, X_{k})$. Since restriction of $x$ on $X_{k-1}$ splits 
along $X_{k}$ by geometric definition of the map 
$\Cal S_{l_{k-1}+1}(X_{k-1})\to LS_{l_k}(F_k)$ we obtain that $\sigma(y)=0$.
Now it follows from commutativity of (3.34)  that $\Theta_{i}(x)=0$.

We prove now the reverse  implication.  
Let $\Theta_{k}(x)=0$. It follows 
from diagram (3.34) that $\Theta_{k-1}(x)=0$. 
Hence there exists an element  
$y\in \Cal S_{n+1}^{k-1}(\Cal X)$ which maps to $x$. 
The last set is identified  
with the classes of concordance of $s$-triangulations 
of the filtartion $\Cal X_{k-1}$. Hence the representative $y$ gives 
an $s$-triangulation    
$(f^{\prime},b^{\prime})$ of $\Cal X_{k-1}$. Since $\Theta_k(x)=0$ 
it follows by the  commutativity (3.34) that $y$ lies in the image
of the map 
$$
\Cal S^k_{n+1}(\Cal X)\to\Cal S_{n+1}^{k-1}(\Cal X)
$$
from (3.34).
Hence $\sigma(x)=0$,  and by decomposition (3.35) the restriction 
of the map $f^{\prime}$ to $X_{k-1}$ splits along
the submanifold $X_k$. We can extend a homotopy to obtain 
an $s$-triangulation of $\Cal X_{k-1}$ for which the restrictriction 
on $X_{k-1}$ is an $s$-triangulation of the pair $(X_k, X_{k-1})$.
Now  application of Proposition 2.5 finishes the proof of the theorem. 
\qed
 \enddemo 
 \bigskip
 
For filtration (1.11) we now 
 describe the relations between  surgery 
obstruction groups 
$LM^i_*(\Cal X)=LM^i(\Cal X_i)$ introduced above 
and stratified $L$-groups $L^{BQ}_*(\Cal X_i)$ 
of Browder and Quinn (see \cite{3} and \cite{35}).     

The Browder-Quinn groups of filtration $\Cal X$ 
 are realized on the spectra level and we recall here 
an inductive definition of these groups from \cite{35, p. 129}
using our notations.  
In accordance with Theorem 2.9  we have a homotopy equivalence
of spectra
$$
\Bbb LM^2(\Cal X)\simeq \Sigma^{l_2-n}\Bbb L^{BQ}(\Cal X_2).
$$ 
It is necessary to remark here that in a similar way a homotopy
equivalence of spectra immediately follows from (2.10), (2.12),
and definition \cite{3, page 129} 
$$
\Bbb LM^1(\Cal X)=\Bbb LP(F_1)\simeq \Sigma^{l_1-n}\Bbb L^{BQ}(\Cal X_1).
\tag 3.36
$$

By Remark 2.3,   filtration $\Cal X$ gives filtration 
$\overline{\Cal X}$ of manifolds  with boundaries. The boundaries 
of the last filtration give a filtration by closed manifolds 
$$
\partial(X_{k-1}\setminus X_k  \subset 
\partial(X_{k-2}\setminus X_k)\subset 
\cdots \subset \partial(X_1 \setminus  X_k)
\subset \partial(X_0\setminus X_k)
\tag 3.37
$$
which we shall denote by $\partial \overline{\Cal X}$. 
Note  that filtrations $\partial \overline{\Cal X}$ and
$\overline{\Cal X}$ contain $k$ spaces, and filtration $\Cal X$ 
contains $k+1$ spaces.

Consider a homotopy commutative diagram of spectra 
$$
\matrix
\Bbb LM^{j}(\overline{\Cal X})& \to & \Bbb LM^{j-1}(\overline{\Cal Y}) &\to &
\Sigma^{l_{j}-n+1} \Bbb L(C)&\to&
 \Sigma \Bbb LM^j(\overline{\Cal X})\\
\downarrow && \downarrow &&\downarrow = && \downarrow\\
\Bbb LM^{j}(\Cal X)& \to & \Bbb LM^{j-1}(\Cal Y) &\to &
\Sigma^{l_{j}-n+1} \Bbb L(C)&\to&
 \Sigma \Bbb LM^j(\Cal X)\\
 \endmatrix
\tag 3.38
$$
where  $j=k-1\geq 1$,
$C=\pi_1(X_0\setminus X_1)$, and horizontal rows are cofibrations by (3.14).
The vertical maps in (3.38) are induced by a natural inclusion 
of filtrations 
$\overline{\Cal X_k}\subset \Cal X_{k-1}$.
For $k-1=j=1$ the central square of (3.38) fits in a homotopy 
commutative diagram of spectra 
$$
\matrix
\Bbb LM^{0}(\overline{\Cal Y}) &\to &
\Sigma^{l_{1}-n+1}\Bbb L(C)\\
\downarrow && \downarrow =\\
\Sigma^{l_1-l_2}\Bbb LM^{1}(\Cal Y) &\to &
\Sigma^{l_{1}-n+1}\Bbb L(C)\\
\downarrow && \downarrow =\\
\Bbb LM^{0}(\Cal Y) &\to &
\Sigma^{l_{1}-n+1} \Bbb L(C)\\
 \endmatrix
\tag 3.39
$$
which follows from diagram (3.14) and diagram (1.8) for the pair
$(X_1, X_2)$ on spectra level.
\smallskip

\proclaim{Proposition 3.10} Let $\Cal X$ be filtration (1.11) with 
$k\geq 2$. There exists the following homotopy
commutative diagram of spectra 
$$
\matrix
\Bbb LM^{k-2}(\overline{\Cal Y}) &\to &
\Sigma^{l_{k-1}-n+1}\Bbb L(C)&\to & \Sigma \Bbb LM^{k-1}(\overline{\Cal X)}\\
\downarrow && \downarrow =&& \downarrow\\
\Sigma^{q_k}\Bbb LM^{k-1}(\Cal Y)&\to&\Sigma^{l_{k-1}-n+1}\Bbb L(C)
&\to& \Sigma^{q_k+1}\Bbb LM^k(\Cal X)\\
\downarrow && \downarrow =&& \downarrow\\
\Bbb LM^{k-2}(\Cal Y) &\to &
\Sigma^{l_{k-1}-n+1} \Bbb L(C)&\to & \Sigma \Bbb LM^{k-1}(\Cal X)\\
 \endmatrix
\tag 3.40
$$
where $l_{k-1}-l_k=q_k$ and $C=\pi_1(X_0\setminus X_1)$. The right 
vertical composition coincides with the right
vertical map in diagram (3.38) for $j=k-1$.
\endproclaim
\demo{Proof} For $k=2$ the result follows from diagrams (3.39) and (3.38)
if we define  
the  right vertical maps in (3.40) as natural maps of homotopical 
 cofibers of the horizontal maps from (3.39) (see \cite{33}). 
Induction on $k$ now
 finishes  proof of the proposition.
 \qed
 \enddemo
 \smallskip
 
\proclaim{Corollary 3.11} Let $\Cal X$ be filtration (1.11) with 
$k\geq 1$. Then a  homotopical 
fiber of the map 
$$
\Bbb LM^{k-1}(\overline{\Cal X})\to 
\Sigma^{q_k}\Bbb LM^k(\Cal X)
\tag 3.41
$$ 
is naturally homotopy equivalent to  $\Sigma^{q_k-1}\Bbb L(\pi_1(X_k))$.
\endproclaim 
\demo{Proof} For $k=1$ the result follows from definition of spectra 
$\Bbb LM^0$, $\Bbb LM^1=\Bbb LP(F_1)$, and cofibration (2.12). For 
$k\geq 2$ the result follows inductively from a homotopy commutative 
diagram 
$$
\matrix
\Sigma^{l_{k-1}-n}\Bbb L(C)&\to & \Bbb LM^{k-1}(\overline{\Cal X)}
&\to&
\Bbb LM^{k-2}(\overline{\Cal Y})  \\
\downarrow= && \downarrow && \downarrow\\
\Sigma^{l_{k-1}-n}\Bbb L(C)
&\to& \Sigma^{q_k}\Bbb LM^k(\Cal X)&\to&
\Sigma^{q_k}\Bbb LM^{k-1}(\Cal Y)\\
\endmatrix
\tag 3.42
$$ 
which follows from diagram (3.40). The right square in (3.42) is pullback
since fibers of horizontal maps are naturally homotopy equivalent. 
Hence fibers of the vertical maps are naturally 
homotopy equivalent, too. By inductive hypothesis fiber 
of the right vertical map is naturally homotopy equivalent 
to $\Sigma^{q_k-1}\Bbb L(\pi_1(X_k))$. From this result of the
corollary follows.\qed
\enddemo
\smallskip

Recall now, that in \cite{35, page 129} an inductive 
definition
of spectra $\Bbb L^{BQ}(\Cal X)$ is given with homotopy groups
$$
\pi_n(\Bbb L^{BQ}(\Cal X))=L_n^{BQ}(\Cal X)
$$
which are Browder-Quinn stratified $L$-groups \cite{3}.

\proclaim{Theorem 3.12} Let $\Cal X$ be the filtration (1.11)
with the smallest manifold of filtration $X_k$ of dimension $l_k$.
We have a naturally homotopy equivalence 
$$
\Bbb LM^k(\Cal X)\simeq \Sigma^{l_k-n}\Bbb L^{BQ}(\Cal X) 
\tag 3.43
$$
 \endproclaim
\demo{Proof} For $k\geq 1$ the cofibration  (3.41) yields a cofibration 
$$
\Bbb L(\pi_1(X_k))\to \Sigma^{1-q_k}\Bbb LM^{k-1}(\overline{\Cal X})\to 
\Sigma^{1}\Bbb LM^k(\Cal X)
$$
which coincides with cofibration for inductive definition 
of the spectra  $\Bbb L^{BQ}(\Cal X)$ \cite{35, page 129} 
up to a shift of dimension of spectra.
\qed
\enddemo
\smallskip

\proclaim{Proposition 3.13} Let $\Cal X$ be filtration (1.11) 
with $k\geq 2$. We have the following braid of exact sequences
$$ 
\matrix 
\rightarrow & L_{l_k+1}(X_k) & 
\rightarrow & 
LM^{k-2}_{l_{k-1}}(\overline{\Cal Y}) & 
         \rightarrow & L_{n-1}(C)& \rightarrow \cr 
    \ & \nearrow \ \ \ \ \ \ \ \ \searrow & \ & \nearrow \ \ \ \ \ 
    \       \searrow 
 & \ & \nearrow \ \ \ \ \ \ \ \ \searrow & \ \cr 
\ & \ & LM_{l_{k-1}}^{k-1}(\overline{\Cal X})& \ & LM^{k-1}_{l_k}(\Cal Y) 
& \ & \ \cr 
                  \ & \searrow \ \ \ \ \ \ \ \ \nearrow & \ & 
                  \searrow \ \ \ \ \ 
                  \ \ 
                  \ 
                  \nearrow 
                  & \ & \searrow \ \ \ \ \ \ \ \ \nearrow & \ \cr 
       \rightarrow & L_{n}(C) & \rightarrow & 
        LM_{l_k}^k(\Cal X) & 
 \rightarrow & L_{l_k}(X_k) & \rightarrow  
\endmatrix 
\tag 3.44
$$ 
where $C=\pi_1(X_0\setminus X_1)$.
\endproclaim 
\demo{Proof} The homotopy long exact sequences of the maps from the 
right pullback 
square of diagram (3.42) give commutative braid of exact 
sequences (3.44). \qed
\enddemo
\bigskip

\subhead 4. Application to Browder-Livesay invariants
\endsubhead
\bigskip

We shall call a  filtration $\Cal X$ (1.11) a  Browder-Livesay filtration 
if for every $k\geq j\geq 1$ the pair of manifolds
$X_{j}\subset X_{j-1}$ is a Browder-Livesay pair. 
Note that $F_j, \ 1\leq j\leq k$,   
is a square of fundamental groups for splitting problem
for the pair of manifolds $X_{j}\subset X_{j-1}$ from  filtration (1.11).
For a Browder-Livesay filtration  any manifold $X_j$ 
is an one-sided submanifold of codimension 1 
in $X_{j-1}$,  horizontal maps in  squares $F_j$ are isomorphisms, 
and vertical maps are inclusions of index 2. 
\smallskip

\proclaim{Theorem 4.1}  Let $\Cal X$ be a Browder-Livesay 
filtration (1.11) in which all the squares $F_j$ for 
$1\leq j\leq k$ are the same. 
Then  filtration of spectra (1.14)  
has the following form 
$$
\Sigma^{k} \Bbb LM^k\to
\Sigma^{k-1} \Bbb LM^{k-1}\to \cdots
\to\Sigma^{1} \Bbb LM^1 \to \Bbb LM^{0}
\tag 4.1
$$
and coincides with the left 
part beginning with $\Bbb X_{0,0}$ of filtration (1.10) 
for spectral sequence of Hambleton and 
Kharshiladze. 
\endproclaim
\demo{Proof} It follows from Corollary (3.2) that the spectrum
$\Bbb LM^j(\Cal X)$ is defined inductively. This spectrum 
is constructed using  three others spectra from 
diagram (3.17) to obtain a pullback square. But the spectrum 
$\Bbb X_{j,0}$ of filtration (1.10) is defined inductively using the    
same construction (see \cite{13} and \cite{26}).\qed
\enddemo
\smallskip

Let $i_+:A\to B^+$ be an inclusion of groups of index 2 as in square 
(1.9). For such inclusion an algebraic version of diagram (1.8) 
was constructed by Ranicki in \cite{31}. It has the following 
form 
$$ 
 \matrix \rightarrow & {L}_{n+1}(A) & \longrightarrow & 
   L_{n+1}(B^+) & 
\overset{\partial}\to{\rightarrow} & LN_{n-1}(A\to B^+)&\rightarrow 
 \cr 
\ & \nearrow \  \  \ \ \ \searrow &\ &^s\nearrow \ \ \ \ \ 
                           \searrow 
   & \ & \nearrow \ \ \  \ \ \searrow & \ \cr 
  \ & \ & L_{n+1}(i^*_-)& \ & L_{n+1}(A\to B^+) & \ & \ \cr 
   \ & \searrow \ \ \ \ \ \nearrow &\ &\searrow \ \ \ \ \ 
\nearrow 
 & \ & \searrow \ \    \ \ \ \nearrow & \ \cr 
\rightarrow & LN_{n}(A\to B^+) & \longrightarrow & 
 L_{n}(B^+) & 
 \longrightarrow & {L}_{n}(A) & \rightarrow 
\endmatrix 
\tag 4.2 
 $$ 
For the Browder-Livesay pair $Y\subset  X$ with the 
square (1.9) of fundamental groups diagram (1.8) coincides 
with diagram  (4.2). The map $\partial $ is called Browder-Livesay
invariant.   If $\partial (x)\neq 0$ then no element 
$x\in L_{n+1}(B)$ can be realized by a normal map of closed 
manifolds (see \cite{9}).   

This diagram is realized on the spectra 
level and we can write down the following pullback square of spectra
$$ 
\matrix 
        &         & \Bbb L(B) &          & \\   
             &\nearrow &           &\searrow  &  \\
\Bbb L(i^*)&         &            &         & \Bbb L(A\to B)   \\
             &\searrow &           &\nearrow  &  \\  
             &         & \Sigma \Bbb L(B^{\epsilon})& & \\
\endmatrix 
\tag 4.3
$$                       
where $i^*$ denotes $i^*_-$ and $B^{\epsilon}$ means that orientation 
on the bottom group $B$ differs from orientation of the upper group 
$B$ outside of the image of the map $i:A\to B$.
Consider a sequence $\Cal A$ of inclusions of subgroups of index 2 into 
a group $B$ with 
an orientation
$$
i_1: A_1\to B; \ i_2:A_2\to B; \ \dots  ; \ i_k: A_k\to B; \ \dots 
\tag 4.4
$$
Every inclusion in (4.4) gives a pullback square similar 
to 4.3 and we can write down the following column of
pullback squares 
$$ 
\matrix 
        &         & \Bbb L(B) &          & \\   
             &\nearrow &           &\searrow  &  \\
\Bbb L(i^*_1)&         &            &         & \Bbb L(A_1\to B)   \\
             &\searrow &           &\nearrow  &  \\  
             &         & \Sigma \Bbb L(B^{\epsilon_1})& & \\
                          &\nearrow &           &\searrow  &  \\
\Sigma\Bbb L(i^*_2)&         &            &         & \Sigma\Bbb L(A_2\to B)\\
             &\searrow &           &\nearrow  &  \\  
             &         & \Sigma^2 \Bbb L(B^{\epsilon_2})& & \\
                          &\nearrow &           &\searrow  &  \\
\Sigma^2\Bbb L(i^*_3)&         &            &         & \Bbb L(A_3\to B)   \\
             &\searrow &           &\nearrow  &  \\  
             &         & \Sigma^3 \Bbb L(B^{\epsilon_3})& & \\
             &         & \vdots              && \\
        &         & \Sigma^{k-1}\Bbb L(B) &          & \\   
             &\nearrow &           &\searrow  &  \\
\Sigma^{k-1}\Bbb L(i^*_k)&   &    &         & \Sigma^{k-1}\Bbb L(A_k\to B) \\
             &\searrow &           &\nearrow  &  \\  
             &         & \Sigma^k \Bbb L(B^{\epsilon_k})& & \\         
                          &         & \vdots              && \\
\endmatrix 
\tag 4.5
$$
in which we have the same agreement on orientations as in
the square (4.3).

Let $\Bbb LM^0(\Cal A)=\Bbb L(B)$ and $\Sigma \Bbb LM^1(\Cal A)=\Bbb L(i_1^*)$.
Using the pullback construction we can extend diagram to the left direction 
similarly to \cite{13}.  In particular, we obtain a  
filtration  of spectra  
$$
\cdots \to \Sigma^{k} \Bbb LM^k(\Cal A) \to
\Sigma^{k-1} \Bbb LM^{k-1}(\Cal A)\to \cdots
\to\Sigma^{1} \Bbb LM^1 (\Cal A)\to \Bbb LM^{0} (\Cal A)
\tag 4.6
$$
which is the left upper diagonal row of the diagram.

We can use filtration (4.6) to construct a surgery spectral 
sequence 
$$
E_r^{p,q}=E_r^{p,q}(\Cal A)
$$ 
for a sequence of inclusions $\Cal A$ 
similar to \cite{13}. We define the first term 
$$
E^{p,q}=\pi_{q-p}(\Sigma^p\Bbb LM^p(\Cal A), \Sigma^{p+1}\Bbb LM^{p+1}(\Cal A)
\cong LN_{q-2p-2}(A_p\to B^{\epsilon_p})  
\tag 4.7
$$
and the first differential 
$$
d_1^{p,q}: E_1^{p,q}\to E_1^{p+1,q} 
$$
which coincides with the composition 
$$
LN_{q-2p-2}(A_{p+1}\to B^{\epsilon_p})\to L_{q-2p-2}(B^{\epsilon_{p+1}})\to  
LN_{q-2p}(A_{p+2}\to B^{\epsilon_{p+1}}).
\tag 4.8
$$  
The first map of the composition 4.8 lies in the diagram (4.3) 
for the inclusion 
$$
A_{p+1}\to B^{\epsilon_{p}},
$$ 
and the second map 
lies in the same diagram for the inclusion 
$$
A_{p+2}\to B^{\epsilon_{p+1}}
$$ 
(see \cite{13}).
Note that the obtained spectral sequence is a natural 
generalization of the spectral 
sequence constructed in \cite{13}. General results about surgery 
spectral sequence from \cite{13} are applicable to obtained 
surgery spectral sequence.
\smallskip

\proclaim{Remark 4.2} Let in (4.4) all the groups $A_i$ be equal to $A$. 
Then the surgery spectral 
sequence constructed above coincides with the 
spectral sequence from  \cite{13} for the inclusion 
$A\to B$
of index
2. 
\endproclaim
\smallskip

Note that a finite sequence $\Cal A$ 
of inclusions $i_j, \  1\leq j\leq k$
as (4.4) gives a finite filtration 
$$
\Sigma^{k} \Bbb LM^k(\Cal A) \to
\Sigma^{k-1} \Bbb LM^{k-1}(\Cal A)\to \cdots
\to\Sigma^{1} \Bbb LM^1 (\Cal A)\to \Bbb LM^{0} (\Cal A)
\tag 4.9
$$
of spectra.
Browder-Livesay filtration $\Cal X$  (1.11)
gives a finite sequence of squares 
$$
F_j=
\left(\matrix 
A_j\ \ \ &\overset{\cong}\to{\to} &A_j\ \ \ \\ 
\downarrow i_- & & \downarrow i_+\\ 
B^{\epsilon_j} &\overset{\cong}\to{\to} &B^{\epsilon_{j-1}}\\ 
\endmatrix\right) 
\tag 4.10
$$ 
of fundamental groups for $ 1\leq j\leq k$ which are similar to (1.9)
Right vertical inclusions from (4.10) give 
finite sequence  $\Cal A(\Cal X)$ of inclusions of index 2 into the
group $B$. 
\smallskip

\proclaim{Proposition 4.3} Under the assumptions above
we have 
$$
\Bbb LM^j(\Cal X) =\Bbb LM^j(\Cal A(\Cal X))
$$
for $1\leq j\leq k$.
\endproclaim
\demo{Proof} The same as the proof of 
Theorem 4.1. \qed 
\enddemo
\bigskip

For the sequence of inclusions  $\Cal A $ (4.4) we can construct 
filtration of spectra (4.6). We denote the homotopy groups
of the spectra  from this filtration as follows
$$
\pi_n(\Bbb LM^k(\Cal A))=LM_n^k(\Cal A).
$$
Thus filtration (4.6) 
gives a tower of groups
$$
\cdots \to LM_{n-j}^j(\Cal A)\to  LM_{n-j+1}^{j-1}(\Cal A)\to \  \cdots \ \to
LM_{n-1}^1(\Cal A)\to LM^0_n(\Cal A)=L_n(B)  
\tag 4.11
$$
Denote by $\phi_j$  the map 
$$
LM_{n-j}^j(\Cal A)\to  LM^0_n(\Cal A)=L_n(B)  
\tag 4.12
$$
given by a composition of maps from (4.11). 
The map $\phi_1$ is the map $s$ in 
diagram (4.2).
\bigskip

\proclaim{Theorem 4.4} Suppose that 
 an element $x\in L_n(B)$,  where $n$ is given by 
$\bmod$ 4, does not  
 lie 
in the the image of $\phi_j$ for some sequence of inclusions 
$\Cal A$ and some natural number $j$.  Then $x$ 
cannot be realized by a normal map of closed manifolds.
\endproclaim
\demo{Proof} Let the element $x\in L_n(B)$ be  realized by 
a normal map of closed manifolds $(f, b):M^n\to X^n$. In accordance 
with \cite{34, \S 9} we can take a product of this surgery 
problem with the projective complex space $P_2(\Bbb C)$ of dimension 4 
to obtain the surgery problem 
$$
f\times Id:M^n\times P_2(\Bbb C) \to X^n\times P_2(\Bbb C)
$$
in dimension $n+4$ with surgery obstruction $x\in L_n(B)$. 
Iterating this construction we can obtain a normal 
map of closed manifolds 
$$
g=f\times Id:M\times (P_2(\Bbb C))^k \to X\times (P_2(\Bbb C))^k=X_0
$$
in dimension $m=n+4k\geq j+5$ with surgery obstruction  
$\Theta_0(g, b^{\prime})=x\in L_n(B)$.
Denote by $\Cal A$  the sequence of inclusions  
$$
i_1: A_1\to B; \ i_2:A_2\to B; \ \dots  ; \ i_j: A_j\to B, 
\tag 4.13
$$
which  defines the map $\phi_j$. Consider a map 
$$
\psi_1:X_0 \to  P_N(\Bbb R), 
$$ 
which  induces epimorphism of fundamental group 
with kernel $A_1$. Here $P_N(\Bbb R)$  is a real projective
space of high dimension. By changing the map $\psi_1$ in 
its homotopy class we can assume that $\psi_1$   
is transversal to $P_{N-1}(\Bbb R)\subset P_N(\Bbb R)$ 
with $\psi_1^{-1}(P_{N-1}(\Bbb R))=X_1$ and that $X_1\subset X_0$ 
is a Browder-Livesay pair (see \cite{5}, \cite{10}, \cite{13}, 
and \cite{17}). Now
in a similar way we consider the map 
$$
\psi_2:X_1 \to  P_N(\Bbb R), 
$$ 
which  induces epimorphism of the fundamental groups 
with kernel $A_2$ with $\psi_2^{-1}(P_{N-1}(\Bbb R))=X_2$ 
and with $X_2\subset X_1$ 
is a Browder-Livesay pair. Iterating this construction 
we obtain a Browder-Livesay filtration $\Cal X$
$$
X_j\subset X_{j-1}\subset\ \cdots \ \subset X_1\subset X_0
\tag 4.14
$$   
with $\Cal A(\Cal X)= \Cal A$. 
From Corollary 3.8 we obtain the following commutative diagram 
$$
\matrix 
H_m(X_0;\bold L_{\bullet}) &\overset{\Theta_0}\to{\longrightarrow}& L_m(B)\\
\ \ \ \ \  \searrow\Theta_j& & \phi_j\nearrow \ \ \ \  \\
&               LM^j_{n-j} & \\
 \endmatrix
 \tag 4.15
 $$
It follows from (4.15) that 
$$
\phi_j\Theta_j(g, b^{\prime})=\Theta_0(g, b^{\prime})=x\in L_m(B)
$$
and hence the element $x$ lies in the image of $\phi_j$.
We have obtained a contradiction and thus the theorem is proved.  
\qed
\enddemo 
\bigskip

An element of $L_n(B)$ doesn't lie in the image $\phi_1$ if and only if 
it maps nontrivially by Browder-Livesay invariant $\partial$
from (4.2). An element  
of $L_n(B)$ doesn't lie in the image $\phi_2$ if and only if 
first or second invariant of Browder-Livesay are nontrivial. 
The second Browder-Livesay invariant was introduced by Hambleton
in \cite{10} (see also \cite{27}) and it is defined only 
if Browder-Livesay invariant is trivial. The iterated Browder-Livesay 
invariants were introduced by Kharshiladze (see \cite{16}, \cite{17},
and \cite{22}). 
Elements of $L_n(B)$ which don't lie in the image 
of $\phi_j$ for some $j$ (and only these elements) 
 detected by iterated Browder-Livesay 
invariants as follows immediately from \cite{13}, \cite{16}, and \cite{17}.
The nonrealizability  of such elements by a normal map
of closed manifolds was proved by Kharshiladze (see  \cite{16} and \cite{17)}
by geometric methods.

\newpage

Authors' addresses:

Yuri V.  Muranov

Department of Informatics and Management

Vitebsk Institute of Modern Knowledge

ul. Gor'kogo 42,

210004  Vitebsk

Belarus

e-mail:  ymuranov\@mail.ru; ymuranov\@imk.edu.by

\bigskip

Rolando Jimenez

Instituto de Matematicas, UNAM

Avenida Universidad S/N, Col. Lomas de Chamilpa

62210 Cuernavaca, Morelos

Mexico

e-mail: rolando\@aluxe.matcuer.unam.mx

\bigskip
Du\v san Repov\v s

Institute for Mathematics, Physics and Mechanics, 
University of Ljubljana,

Jadranska 19, Ljubljana, 

Slovenia

e-mail: dusan.repovs\@uni-lj.si
\enddocument
\bye